\documentclass[10pt]{article}
\usepackage{amssymb,latexsym}
\usepackage{epsfig}
\usepackage{eufrak}
\usepackage{amsmath}
\usepackage{mathrsfs}
\usepackage{color}

\newtheorem{theorem}{Theorem}
\newtheorem{lemma}{Lemma}

\newcommand{\be}{\begin{equation}}
\newcommand{\ee}{\end{equation}}
\newcommand{\bee}{\begin{eqnarray*}}
\newcommand{\eee}{\end{eqnarray*}}
\newcommand{\bel}{\begin{eqnarray}}
\newcommand{\eel}{\end{eqnarray}}
\newcommand{\bec}{\begin{cases}}
\newcommand{\eec}{\end{cases}}
\newcommand{\bem}{\begin{bmatrix}}
\newcommand{\eem}{\end{bmatrix}}

\newcommand{\la}{\label}
\newcommand{\li}{\left}
\newcommand{\ri}{\right}

\newcommand{\ovl}{\overline}

\newcommand{\lm}{\lambda}
\newcommand{\Lm}{\Lambda}
\newcommand{\Up}{\Upsilon}

\newcommand{\si}{\sigma}
\newcommand{\Si}{\Sigma}

\newcommand{\vDe}{\varDelta}

\newcommand{\ga}{\gamma}
\newcommand{\Ga}{\Gamma}
\newcommand{\vse}{\vartheta}
\newcommand{\se}{\theta}
\newcommand{\Se}{\Theta}

\newcommand{\al}{\alpha}
\newcommand{\ba}{\beta}

\newcommand{\ro}{\rho}

\newcommand{\om}{\omega}
\newcommand{\Om}{\Omega}

\newcommand{\f}{\frac}
\newcommand{\sq}{\sqrt}
\newcommand{\cd}{\cdots}

\newcommand{\qu}{\quad}
\newcommand{\qqu}{\qquad}
\newcommand{\fa}{\forall}

\newcommand{\mscr}{\mathscr}
\newcommand{\mcal}{\mathcal}

\newcommand{\bb}{\mathbb}
\newcommand{\fra}{\mathfrak}

\newcommand{\wh}{\widehat}

\newcommand{\bs}{\boldsymbol}

\newcommand{\sh}{\slash}

\newcommand{\tx}{\text}

\newcommand{\iy}{\infty}

\newcommand{\bed}{\begin{description}}
\newcommand{\eed}{\end{description}}
\newcommand{\bei}{\begin{itemize}}
\newcommand{\eei}{\end{itemize}}
\newcommand{\ben}{\begin{enumerate}}
\newcommand{\een}{\end{enumerate}}
\newcommand{\bib}{\bibitem}
\newcommand{\beL}{\begin{lemma}}
\newcommand{\eeL}{\end{lemma}}
\newcommand{\beT}{\begin{theorem}}
\newcommand{\eeT}{\end{theorem}}
\newcommand{\sect}{\section}

\newcommand{\bpf}{\begin{pf}}
\newcommand{\epf}{\end{pf}}

\newcommand{\bi}{\binom}

\newcommand{\pfbox}{\hfill\mbox{$\Box$}}

\newenvironment{pf}{\paragraph*{Proof{\rm.}}}{\pfbox\bigskip}

\begin{document}

\title{{\bf Concentration Inequalities from Likelihood Ratio Method}
\thanks{The author is afflicted  with the Department
of Electrical Engineering and Computer Science at Louisiana State University, Baton Rouge, LA 70803, USA, and the Department of Electrical
Engineering, Southern University and A\&M College, Baton Rouge, LA 70813, USA; Email: xinjiachen@lsu.edu}}

\author{Xinjia Chen}

\date{September  2014}

\maketitle

\begin{abstract}

We explore the applications of  our previously established likelihood-ratio method for deriving concentration inequalities for a wide variety of
univariate and multivariate distributions. New concentration inequalities for various distributions are developed without the idea of minimizing
moment generating functions.

\end{abstract}

\tableofcontents

\sect{Introduction}

Bounds for probabilities of random events play important roles in many areas of engineering and sciences.  Formally, let $E$ be an event defined
in probability space $(\Om, \Pr , \mscr{F})$, where $\Om$ is the sample space, $\Pr$ denotes the probability measure, and $\mscr{F}$ is the
$\si$-algebra.  A frequent problem is to obtain simple bounds as tight as possible for $\Pr \{ E \}$. In general, the event $E$ can be expressed
in terms of a matrix-valued random variable $\bs{\mcal{X}}$.  In particular, $\bs{\mcal{X}}$ can be a random vector or scalar. Clearly, the
event $E$ can be represented as $\{ \bs{\mcal{X}} \in \mscr{E} \}$, where $\mscr{E}$ is a certain set of deterministic matrices. In probability
theory, a conventional approach for deriving inequalities for  $\Pr \{ E \}$ is to bound the indicator function $\bb{I}_{ \{ \bs{\mcal{X}} \in
\mscr{E} \} }$ by a family of random variables having finite expectation and minimize the expectation.  The central idea of this approach is to
seek a family of bounding functions $w(\bs{\mcal{X}}, \bs{\vse})$ of $\bs{\mcal{X}}$, parameterized by $\bs{\vse} \in \varTheta$, such that \be
\la{inviewvip}
 \bb{I}_{ \{  \bs{\mcal{X}} \in \mscr{E}   \} } \leq
w(\bs{\mcal{X}}, \bs{\vse}) \qu \tx{for all $\bs{\vse} \in
\varTheta$}. \ee Here, the notion of inequality (\ref{inviewvip}) is
that the inequality $\bb{I}_{ \{  \bs{\mcal{X}} (\om) \in \mscr{E}
\} } \leq w(\bs{\mcal{X}} (\om) , \bs{\vse})$ holds for every $\om
\in \Om$. As a consequence of the monotonicity of the mathematical
expectation $\bb{E} [. ]$, \be \la{inviewvip889b}
 \Pr \{  E \}
= \bb{E} [ \bb{I}_{ \{  \bs{\mcal{X}} \in \mscr{E}   \} }  ] \leq \bb{E} [ w(\bs{\mcal{X}}, \bs{\vse}) ] \qu \tx{for all $\bs{\vse} \in
\varTheta$}. \ee Minimizing the upper bound in (\ref{inviewvip889b}) with respect to $\bs{\vse} \in \varTheta$ yields  \be \la{inviewvip889c}
\Pr \{  E \} \leq \inf_{ \bs{\vse} \in \varTheta } \bb{E} [ w(\bs{\mcal{X}}, \bs{\vse}) ]. \ee Classical concentration inequalities such as
Chebyshev inequality and  Chernoff bounds \cite{Chernoff} can be derived by this approach with various bounding functions $w(\bs{\mcal{X}},
\bs{\vse})$, where $\bs{\mcal{X}}$ is a scalar random variable. We call this technique of deriving probabilistic inequalities as the {\it
mathematical expectation} (ME) method, in view of the crucial role played by the mathematical expectation of bounding functions.  For the ME
method to be successful,  the mathematical expectation $\bb{E} [ w (\bs{\mcal{X}}, \bs{\vse}) ]$ of the family of bounding functions  $w
(\bs{\mcal{X}}, \bs{\vse}), \; \bs{\vse} \in \varTheta$ must be convenient for evaluation and minimization.  The ME method is a very general
approach. However, it has two drawbacks. First, in some situations, the mathematical expectation $\bb{E} [ w (\bs{\mcal{X}}, \bs{\vse}) ]$ may
be intractable. Second, the ME method may not fully exploit the information of the underlying distribution, since the mathematical expectation
is only a quantity of summary for the distribution.

Recently,  we have proposed in \cite{ChenLR, ChenLRgeneral, ChenLRSPIR}  a more general approach for deriving probabilistic inequalities, aiming
at overcoming the drawbacks of the ME method. Let $f(.)$ denote the probability density function (pdf) or probability mass function (pmf) of
$\bs{\mcal{X}}$. The primary idea of the proposed approach is to seek a family of pdf or pmf $g(., \bs{\vse})$, parameterized by $\bs{\vse} \in
\varTheta$,   and a deterministic function $\Lm(\bs{\vse})$ of $\bs{\vse} \in \varTheta$  such that  for all $\bs{\vse} \in \varTheta$,  the
indicator function $\bb{I}_{ \{ \bs{\mcal{X}} \in \mscr{E} \} }$ is bounded from above by the product of $\Lm(\bs{\vse})$ and the likelihood
ratio $\f{ g(\bs{\mcal{X}}, \bs{\vse}) }{ f(\bs{\mcal{X}}) }$. Then, the probability $\Pr \{ \bs{\mcal{X}} \in \mscr{E} \}$ is bounded from
above by the infimum of $\Lm(\bs{\vse})$ with respect to $\bs{\vse} \in \varTheta$.  Due to the central role played by the likelihood ratio,
this technique of deriving probabilistic inequalities is referred to as the {\it likelihood ratio} (LR) method.  It has been demonstrated in
\cite{ChenLRgeneral} that the ME method is actually a special technique of the LR method.

In this paper, we shall apply the LR method to investigate the concentration phenomenon of random variables. Our goal is to derive simple and
tight concentration inequalities for various distributions.  The remainder of the paper is organized as follows.  In Section 2, we introduce the
fundamentals of the LR method. In Section 3, we apply the LR method to the development of concentration inequalities for univariate
distributions. In Section 4, we apply the LR method to establish concentration inequalities for multivariate distributions. Section 5 is the
conclusion. Most proofs are given in Appendices.

Throughout this paper, we shall use the following notations.   Let $\bb{I}_E$ denote the indicator function such that $\bb{I}_E = 1$ if $E$ is
true and $\bb{I}_E = 0$ otherwise.   We use the notation $\bi{t}{k}$ to denote a generalized combinatoric number in the sense that
\[
\bi{t}{k} = \f{\prod_{\ell = 1}^k (t -\ell + 1)}{k!} = \f{ \Ga(t + 1) }{\Ga(k + 1) \; \Ga(t - k + 1)},  \qqu  \bi{t}{0} = 1,
\]
where $t$ is a real number and $k$ is a non-negative integer.  We use $\ovl{X}_n$ to denote the average of random variables $X_1, \cd, X_n$,
that is, $\ovl{X}_n = \f{\sum_{i=1}^n X_i}{n}$.  The notation $\top$ denotes the transpose of a matrix.  The trace of a matrix is denoted by
$\tx{tr}$.  We use pdf and pmf to represent probability density function and probability mass function, respectively.   The other notations will
be made clear as we proceed.

\sect{Likelihood Ratio Method}

In this section, we shall introduce the LR method for deriving probabilistic inequalities.

\subsection{General Principle}

 Let $E$ be an event which can be expressed in terms of matrix-valued
random variable $\bs{\mcal{X}}$, where $\bs{\mcal{X}}$ is defined on the sample space $\Om$ and $\si$-algebra $\mscr{F}$ such that the true
probability measure is one of two measures $\Pr$ and $\bb{P}_{\bs{\vse} }$. Here, the measure $\Pr$ is determined by pdf or pmf $f(.)$. The
measure $\bb{P}_{\bs{\vse} }$ is determined by pdf or pmf $g(., \bs{\vse})$, which is parameterized by $\bs{\vse} \in \varTheta$.  The subscript
in $\bb{P}_{\bs{\vse}}$ is used to indicate the dependence on the parameter $\bs{\vse}$. Clearly, there exists a set, $\mscr{E}$, of
deterministic matrices of the same size as $\bs{\mcal{X}}$ such that $E = \{ \bs{\mcal{X}} \in \mscr{E} \}$. The LR method for obtaining an
upper bound for the probability $\Pr \{ E \}$ is based on the following general result. \beT \la{ThM888} Assume that there exists a function
$\Lm(\bs{\vse})$ of $\bs{\vse} \in \varTheta$ such that \be \la{LRBVIP} f (\bs{\mcal{X}}) \; \bb{I}_{ \{ \bs{\mcal{X}} \in \mscr{E} \} } \leq
\Lm(\bs{\vse}) \; g (\bs{\mcal{X}}, \bs{\vse}) \qu \tx{for all $\bs{\vse} \in \varTheta$}. \ee Then, \be \la{mainVIP}
 \Pr \{  E   \} \leq \inf_{ \bs{\vse} \in \varTheta }  \Lm(\bs{\vse})  \; \bb{P}_{\bs{\vse} } \{ E \}
  \leq \inf_{ \bs{\vse} \in \varTheta } \Lm(\bs{\vse}). \ee  In particular,
  if the infimum of $\Lm (\bs{\vse})$ is attained at $\vse^* \in
  \varTheta$, then
  \be
  \la{maintight}
\Pr \{  E   \} \leq \bb{P}_{\vse^*} \{ E \} \; \Lm(\vse^*).
  \ee

  \eeT

The notion of the inequality in (\ref{LRBVIP}) is that $f (\bs{\mcal{X}} (\om) )  \; \bb{I}_{ \{ \bs{\mcal{X}} (\om)  \in \mscr{E} \} } \leq
\Lm(\bs{\vse}) \; g (\bs{\mcal{X}} (\om), \bs{\vse})$ for every $\om \in \Om$.   The function $\Lm(\bs{\vse})$ in (\ref{LRBVIP}) is referred to
as {\it likelihood-ratio bounding function}. Theorem \ref{ThM888} asserts that the probability of event $E$ is no greater than the likelihood
ratio bounding function.

\subsection{Construction of Parameterized Distributions}

In the sequel, we shall introduce two approaches for constructing parameterized distributions $g(., \bs{\vse})$ which are essential for the
application of the LR method.

\subsubsection{Weight Function} \la{weightpara}

A natural approach to construct parameterized distribution $g(.,
\bs{\vse})$ is to modify the pdf or pmf $f(.)$ by multiplying it
with a parameterized function and performing a normalization.
Specifically, let $w (.,  \bs{\vse} )$ be a non-negative function
with parameter $\bs{\vse} \in \varTheta$ such that $\bb{E} [ w
(\bs{\mcal{X}}, \bs{\se}) ]  < \iy$ for all $\bs{\vse} \in
\varTheta$,  where the expectation is taken under the probability
measure $\Pr$ determined by $f(.)$. Define a family of distributions
as
\[ g(\bs{\mcal{X}}, \bs{\vse}) =  \f{ w (\bs{\mcal{X}}, \bs{\vse}) \; f( \bs{\mcal{X}} )  }{\bb{E} [ w (\bs{\mcal{X}}, \bs{\vse}) ] }
\]
for  $\bs{\vse} \in \varTheta$ and $\bs{\mcal{X}}$ in the range of
$\bs{\mcal{X}}$.  In view of its role in the modification of $f(.)$
as $g(., \bs{\vse})$, the function $w (. , \bs{\vse})$ is called a
{\it weight function}. Note that \be \la{nowgood} f(\bs{\mcal{X}})
\; w (\bs{\mcal{X}}, \bs{\vse})  = \bb{E} [ w (\bs{\mcal{X}},
\bs{\vse}) ] \; g(\bs{\mcal{X}}, \bs{\vse}) \qu \tx{for all
$\bs{\vse} \in \varTheta$}. \ee For simplicity, we choose the weight
function such that the condition (\ref{inviewvip}) is satisfied.
Combining (\ref{inviewvip}) and (\ref{nowgood}) yields  \[
f(\bs{\mcal{X}}) \; \bb{I}_{ \{ \bs{\mcal{X}} \in \mscr{E} \} } \leq
f(\bs{\mcal{X}}) \; w (\bs{\mcal{X}}, \bs{\vse}) = \bb{E} [ w
(\bs{\mcal{X}}, \bs{\vse}) ] \; g(\bs{\mcal{X}}, \bs{\vse}) \qu
\tx{for all $\bs{\vse} \in \varTheta$}. \]   Thus, the likelihood
ratio bounding function can be taken as
\[
\Lm (\bs{\vse}) = \bb{E} [ w (\bs{\mcal{X}}, \bs{\vse})  ] \qu
\tx{for $\bs{\vse} \in \varTheta$}.
\]
It follows from Theorem \ref{ThM888} that  \[ \Pr \{  E  \} \leq
\inf_{ \bs{\vse} \in \varTheta }  \Lm (\bs{\vse}) \; \bb{P}_{\bs{\vse}} \{ E \}
\leq \inf_{ \bs{\vse} \in \varTheta } \Lm (\bs{\vse}).
\]
Thus, we have demonstrated that the ME method  is actually a special technique of the LR method.

By constructing a family of parameterized distributions and making use of the LR method, we have obtained the following result.

\beT

\la{ChenChBE}

Let $X$ be a random variable with moment generating function
$\phi(.)$. Let $X_1, \cd, X_n$ be i.i.d. samples of $X$.  Let
$C_{BE}$ be the absolute constant in the Berry-Essen inequality.
Then,
\[
\Pr \{ \ovl{X}_n \geq z \} \leq \li ( \f{1}{2} + \vDe \ri ) \li [ e^{- z \tau} \phi(\tau) \ri ]^n,
\]
where \[ \vDe = \min \li \{ \f{1}{2}, \;   \f{ C_{BE} }{\sq{n}} \li
( \f{ \phi (\tau) [ \phi^{\prime \prime \prime \prime} (\tau) - 4 z
\phi^{ \prime \prime \prime } (\tau) ] + 3 [ \phi^{ \prime \prime }
(\tau) ]^2 } { [ \phi^{\prime \prime} (\tau) - z^2 \phi(\tau)  ]^2 }
- 3 \ri )^{\f{3}{4}} \ri \}
\]
with $\tau$ satisfying $\f{\phi^\prime(\tau)}{ \phi(\tau) } = z$.

\eeT

See Appendix \ref{ChenChBEapp} for a proof.  Note that $\vDe \to 0$ as $n \to \iy$.
So, for large $n$, the above bound is twice tighter than
the classical Chernoff bound.

\subsubsection{Parameter Restriction}

In many situations, the pdf or pmf $f(.)$ of $\bs{\mcal{X}}$ comes
from a family of distributions parameterized by $\bs{\se} \in \Se$.
If so, then the parameterized distribution $g(., \bs{\vse})$ can be
taken as the subset of pdf or pmf with parameter $\bs{\vse}$
contained in a subset $\varTheta$ of parameter space $\Se$.  By
appropriately choosing the subset $\varTheta$, the deterministic
function $\Lm(\bs{\vse})$ may be readily obtained.  As an
illustrative example, consider the normal distribution.

A random variable $X$ is said to have a normal distribution
with mean $\mu$ and variance $\si^2$  if it possesses a probability density function
\[
    f_X (x) = \frac{1}{\sqrt{2 \pi} \sigma} \exp \li ( -\f{|x-\mu|^2}{2\sigma^2} \ri ).
\]
Let $X_1, \cd, X_n$ be i.i.d. samples of the random variable $X$.
The following well-known inequalities hold true.
 \bel
&  & \Pr \{ \ovl{X}_n \leq z \} \leq \f{1}{2} \exp \li ( - \f{n (z -
\mu)^2}{2 \si^2} \ri ) \qu \tx{for $z \leq \mu$},  \la{norm1}\\
&  & \Pr \{ \ovl{X}_n \geq z \} \leq \f{1}{2} \exp \li ( - \f{n (z -
\mu)^2}{2 \si^2} \ri ) \qu \tx{for $z \geq \mu$}. \la{norm2} \eel
 It should be noted that the factor $\f{1}{2}$ in these inequalities
cannot be obtained by using conventional techniques of Chernoff
bounds. By virtue of the LR method, we can provide an easy proof for
inequalities (\ref{norm1}) and (\ref{norm2}).  We proceed as
follows.

Let $\bs{\mcal{X}} = [ X_1, \cd, X_n ]$ and $\bs{x} = [x_1, \cd,
x_n]$.  The joint
 probability density function of $\bs{\mcal{X}}$ is
 \[
f_{\bs{\mcal{X}}} (\bs{x}) = \frac{1}{ (\sqrt{2 \pi} \sigma)^n }
\exp \li (  -\f{ \sum_{i=1}^n ( x _i -\mu )^2}{2\sigma^2} \ri ).
 \]
To apply the LR method to show (\ref{norm1}), we construct a family of probability density functions
\[
g_{\bs{\mcal{X}}} (\bs{x}, \vse) = \frac{1}{ (\sqrt{2 \pi} \sigma)^n
} \exp \li (  -\f{ \sum_{i=1}^n ( x _i  - \vse)^2}{2\sigma^2} \ri )
 \]
for $\vse \in  (- \iy, z]$ with $z \leq \mu$. It can be checked that
\[
\f{ f_{\bs{\mcal{X}}} (\bs{x}) }{ g_{\bs{\mcal{X}}} (\bs{x}, \vse) }
= \li [ \exp \li (  - \f{ 2 (\vse - \mu)  \ovl{\bs{x}}_n  + \mu^2 -
\vse^2}{2\sigma^2} \ri ) \ri ]^n.
\]
For any $\vse \in (- \iy, z]$, we have $\vse \leq z \leq \mu$ and thus
\[
\f{ f_{\bs{\mcal{X}}} (\bs{x}) }{ g_{\bs{\mcal{X}}} (\bs{x}, \vse) } \leq \li [ \exp \li (  - \f{ 2 (\vse - \mu) z + \mu^2 - \vse^2}{2\sigma^2}
\ri ) \ri ]^n \qu \fa \vse \in (- \iy, z] \; \tx{for $\ovl{\bs{x}}_n \leq z$}.
\]
This implies that
\[
f_{\bs{\mcal{X}}} ( \bs{\mcal{X}}) \; \bb{I}_{ \{ \ovl{X}_n \leq z \} } \leq \Lm (\vse) \; g_{\bs{\mcal{X}}} ( \bs{\mcal{X}}, \vse) \qu \fa \vse
\in (- \iy, z],
\]
where
\[
\Lm (\vse) = \li [ \exp \li (  - \f{ 2 (\vse - \mu) z + \mu^2 -
\vse^2}{2\sigma^2} \ri ) \ri ]^n.
\]
By differentiation, it can be readily shown that the infimum of $\Lm (\vse)$
with respect to $\vse \in (- \iy, z]$ is  equal to \[ \exp \li ( -
\f{n (z - \mu)^2}{2 \si^2} \ri ),
\]  which is attained at $\vse = z$. By symmetry, it can be shown that
\[
\bb{P}_z  \{ \ovl{X}_n \leq z \} = \f{1}{2}.
\]
Using these facts and invoking (\ref{maintight}) of Theorem 1, we
have
\[
\Pr \{ \ovl{X}_n \leq z \} \leq \bb{P}_z  \{ \ovl{X}_n \leq z \} \;
\Lm (z) \qu \tx{for $z \leq \mu$}.
\]
This implies that inequality (\ref{norm1}) holds. In a similar
manner, we can show inequality (\ref{norm2}).

\section{Concentration Inequalities for Univariate Distributions}

In this section, we shall apply the LR method to derive bounds for tail probabilities for univariate distributions. Such bounds are referred to
as concentration inequalities.

\subsection{Beta Distribution}

A random variable $X$ is said to have a beta distribution if it
possesses a probability density function
\[
f (x) = \f{ 1 }{ \mcal{B} (\al, \ba)  } x^{\al - 1} (1 - x)^{\ba - 1}, \qqu 0 < x  < 1, \qu \al > 0, \qu \ba > 0,
\]
where $\mcal{B} (\al, \ba) = \f{ \Ga(\al) \Ga(\ba)  }{ \Ga(\al +
\ba) }$.  Let $X_1, \cd, X_n$ be i.i.d. samples of the random
variable $X$. Making use of the LR method, we have shown the
following results.

 \beT

\la{betathm}

Let $z \in (0, 1)$ and $\mu = \bb{E} [ X ] = \f{\al}{\al + \ba}$.
Define $\wh{\al} = \f{\ba z}{1 - z}$ and $\wh{\ba} = \f{\al (1 - z)}{z}$.
Then,  \bel &  &  \Pr \li \{ \ovl{X}_n \leq z \ri \} \leq \li [  \f{ \mcal{B} (\wh{\al}, \ba) } {
\mcal{B} (\al, \ba) } \f{z^\al}{z^{\wh{\al}}} \ri ]^n \qqu \tx{for $0 < z \leq \mu$}, \la{beta888a}\\
&  & \Pr \li \{ \ovl{X}_n \geq z \ri  \} \leq \li [ \f{ \mcal{B}
(\al, \wh{\ba}) } {  \mcal{B} (\al, \ba) } \f{(1 - z)^\ba}{(1 -
z)^{\wh{\ba}}} \ri ]^n \qqu \tx{for $\mu \leq  z < 1$ }.
\la{beta888b} \eel Specially, if $\ba = 1$, then \be  \Pr \li \{
\ovl{X}_n \leq z \ri  \} \leq  \li ( e \al z^\al \ln \f{1}{z} \ri
)^n \qu \tx{for $0 < z < \exp \li ( - \f{1}{\al} \ri )$}.
\la{beta888c} \ee

\eeT

See Appendix \ref{betathmapp} for a proof.

\subsection{Beta Negative Binomial Distribution}

A random variable $X$ is said to have a beta distribution if it
possesses a probability mass function
\[
f (x) = \Pr \{ X = x \} = \bi{n + x - 1}{x} \f{ \Ga ( \al + n)
\Ga(\ba + x) \Ga (\al + \ba)  }
 {  \Ga (\al + \ba + n + x) \Ga(\al) \Ga(\ba) }, \qqu x =
0, 1, 2, \cd
\]
where $\al > 1$ and $\ba > 0$ and $n > 1$.  By virtue of the LR
method, we have obtained the following results.

\beT

\la{ineqthm4}

Let $z$ be a nonnegative integer no greater than $\f{ n \ba }{\al -
1 }$.   Then,
\[
\Pr \{ X \leq z \} \leq \f{ \Ga( \f{ \al z - z }{n} )  }{ \Ga(\ba) }
\f{ \Ga(\ba + z) } { \Ga(\f{ \al z - z }{n} + z)  } \f{ \Ga (\al +
\f{ \al z - z }{n} + n + z) } { \Ga (\al + \ba + n + z) }.
 \]

\eeT

See Appendix \ref{ineqthm4app} for a proof.

\subsection{Beta-Prime Distribution}

A random variable $X$ is said to have a beta-prime distribution if
it possesses a probability density function
\[
f (x) = \f{ x^{\al - 1} (1 + x)^{-\al - \ba}  }{ \mcal{B} (\al, \ba) }, \qqu x > 0, \qu \al > 0, \qu \ba > 0.
\]
Let $X_1, \cd, X_n$ be i.i.d. samples of the random variable $X$.
Making use of the LR method, we have obtained the following results.

\beT

\la{betaprime}

Assume that $\ba > 1$ and $0 < z \leq \f{\al}{\ba - 1}$. Then, \bel
& &
 \Pr \li \{ \ovl{X}_n \leq z \ri \} \leq
 \li [  \li (\f{z}{1 + z} \ri )^{\al + z - \ba z} \f{ \mcal{B}( \ba z - z, \ba) } { \mcal{B}(\al, \ba) } \ri
 ]^n, \la{betaprimea}\\
&  &  \Pr \li \{ \ovl{X}_n \leq z \ri \} \leq \li [ \f{
\mcal{B}(\al, 1 + \f{\al}{z}) }{ \mcal{B}(\al, \ba) } (1 + z)^{1 +
\f{\al}{z} - \ba} \ri ]^n. \la{betaprimeb} \eel

\eeT

See Appendix \ref{betaprimeapp} for a proof.

\subsection{Borel Distribution}

A random variable $X$ is said to possess a Borel distribution if it
has a probability mass function
\[
f (x) = \Pr \{ X = x \} = \f{(\se x)^{x-1} e^{- \se x}}{x!}, \qqu x = 1, 2, \cd,
\]
where $0 < \se < 1$.  Let $X_1, \cd, X_n$ be i.i.d. samples of the
random variable $X$. Making use of the LR method, we have obtained
the following result.

\beT

\la{borelcdf}

\be \la{borelineq} \Pr \{ \ovl{X}_n \leq z \} \leq \li [ \li (  \f{
e \se z }{ 1 - z } \ri )^{z - 1} e^{- \se z} \ri ]^n \qu \tx{for $1
< z < \f{1}{1 - \se}$}. \ee

\eeT

See Appendix \ref{borelcdfapp} for a proof.

\subsection{Consul Distribution}

A random variable $X$ is said to have a Consul distribution if it
possesses a probability mass function
\[
f (x) = \Pr \{ X = x \} = \f{1}{x} \bi{mx}{x-1} \li ( \f{\se}{1 - \se} \ri )^{x - 1} (1 - \se)^{m x}, \qqu x = 1, 2, \cd
\]
where $0 \leq \se <1, \; 1 \leq m < \f{1}{\se}$. See, e.g., \cite{Consul}, for an introduction of this distribution.   Let $X_1, \cd, X_n$ be
i.i.d. samples of the random variable $X$. Making use of the LR method, we have obtained the following result.

\beT

\la{consulcdf}

\be \la{consineq}
 \Pr \{ \ovl{X}_n \leq z \} \leq  \li [ \f{ \li (
\f{\se}{1 - \se} \ri )^{z - 1} (1 - \se)^{m z}  }
 {  \li (  \f{z - 1}{1 - z + m z} \ri )^{z - 1} (1 - \f{z - 1}{m z})^{m z}  } \ri ]^n \qu \tx{for $1 \leq z < \f{1}{1 - m \se}$}.
\ee

\eeT

See Appendix \ref{consulcdfapp} for a proof.

\subsection{Geeta Distribution}

A random variable $X$ is said to have a Geeta distribution if it
possesses a probability mass function
\[
f (x) = \Pr \{ X = x \} = \f{1}{\ba x - 1} \bi{\ba x - 1}{x}
\se^{x-1} (1 - \se)^{\ba x - x}, \qqu x = 1, 2, \cd
\]
where $0 < \se < 1$ and $1 < \ba < \f{1}{\se}$.  Making use of the LR method, we have obtained the following result.

\beT \la{ineqgeeta}

\be \la{geetaineq}
 \Pr \{ \ovl{X}_n \leq z \} \leq \li [  \f{ \se^{z-1} (1 - \se)^{\ba z - z}  }{ \li ( \f{z -
1}{\ba z - 1} \ri )^{z-1} \li ( 1 - \f{z - 1}{\ba z - 1} \ri )^{\ba
z - z} } \ri ]^n \qu \tx{for} \; 1 \leq z \leq \f{1 - \se}{1 - \ba
\se}.
\ee

\eeT

See Appendix \ref{ineqgeetaapp} for a proof.

\subsection{Gumbel Distribution}

A random variable $X$ is said to have a Gumbel distribution if it
possesses a probability density function
 \[
f (x) = \f{1}{\ba} \exp \li [ \f{\mu - x}{\ba}  - \exp \li ( \f{\mu
- x}{\ba} \ri ) \ri ], \qqu  - \iy < x < \iy,
 \]
where $\ba > 0$ and $- \iy < \mu < \iy$.  Let $X_1, \cd, X_n$ be
i.i.d. samples of random variable $X$. By virtue of the LR method,
we have obtained the following result.

\beT

\la{ineqgumbel}

\be \la{gumbelineq}
  \Pr \{ \ovl{X}_n \leq z \}  \leq  \li \{ \exp
\li [ \f{ \mu - z }{\ba} + 1 - \exp \li ( \f{ \mu - z }{\ba} \ri )
\ri ] \ri \}^n \ee for $z \leq \mu$.

\eeT

See Appendix \ref{ineqgumbelapp} for a proof.

\subsection{Inverse Gamma Distribution}

A random variable $X$ is said to have an inverse gamma distribution
if it possesses a probability density function
\[
f(x) = \f{ \ba^\al } {  \Ga ( \al) } x^{-\al - 1} \exp \li ( -
\f{\ba}{x} \ri ), \qqu x > 0, \qqu  \al > 0, \qu \ba > 0.
\]
Let $X_1, \cd, X_n$ be i.i.d. samples of random variable $X$.  By
virtue of the LR method, we have obtained the following results.

\beT

\la{invGa}

\bel \Pr \{ \ovl{X}_n \leq z \} & \leq &  \li [ \f{ \Ga (\f{\ba}{z}
+ 1) }{ \Ga (\al) }   \li (  \f{ z } { \beta } \ri )^{\f{\ba}{z} -
\al + 1} \ri ]^n  \qu \tx{for $0 < z \leq \f{ \ba } { \al - 1 }$}, \la{invgamma1}\\
\Pr \{ \ovl{X}_n \leq z \} & \leq &  \li [ \li ( \f{\ba}{\al z} \ri
)^\al
 \exp \li ( \f{ \al z - \ba}{ z } \ri ) \ri ]^n \qu \tx{for $0 < z \leq
\f{\ba}{\al}$}. \la{invgamma2} \eel

\eeT

See Appendix \ref{invGaapp} for a proof.

\subsection{Inverse Gaussian Distribution}

A random variable $X$ is said to have an inverse Gaussian
distribution if it possesses a probability density function
\[
f(x) = \li ( \f{\lm}{2 \pi x^3}  \ri )^{1 \sh 2}  \exp \li ( - \f{ \lm (x - \se)^2}{2 \se^2 x}   \ri ), \qqu x > 0
\]
where $\lm > 0$ and $\se > 0$.

Let $X_1, \cd, X_n$ be i.i.d. samples of random variable $X$.  By virtue of the LR method, we have obtained the following result.

\beT

\la{ineqinvGaussian}

\be \la{inversegausineq}
 \Pr \{ \ovl{X}_n \leq z \}  \leq  \li [ \exp \li (
\f{\lm}{\se} - \f{\lm}{2 z} - \f{ \lm z }{ 2 \se^2  }  \ri ) \ri ]^n
\qu \tx{for} \; 0 < z \leq \se. \ee

\eeT

See Appendix \ref{ineqinvGaussianapp} for a proof.

\subsection{Lagrangian Logarithmic Distribution}

A random variable $X$ is said to have a Lagrangian logarithmic
distribution if it possesses a probability mass function
\[
f(x) = \Pr \{ X = x\} = \f{ - \se^x (1 - \se)^{x (\ba - 1)} \Ga (\ba x)  }
 { \Ga(x + 1) \Ga (\ba x - x + 1) \ln (1 - \se) }, \qqu x = 1, 2, \cd
\]
where $0 < \se \leq \se \ba < 1$.  Let $X_1, \cd, X_n$ be i.i.d.
samples of random variable $X$.  By virtue of the LR method, we have obtained
the following result.

\beT

\la{Laglog}

\be \la{LagBin}
 \Pr \{  \ovl{X}_n \leq z \} \leq \li [ \li (  \f{\se}{\vse}  \ri )^z  \li (  \f{1 - \se}{1 - \vse}  \ri
 )^{z (\ba - 1)} \f{ \ln ( 1 - \vse) } { \ln ( 1 - \se) }
\ri ]^n \qu \tx{for} \; 0 < z \leq \f{ \se }{ (\ba \se - 1) \ln ( 1
- \se) }, \ee
 where $\vse$
satisfies the equation $z = \f{ \vse }{ (\ba \vse - 1) \ln ( 1 -
\vse) }$.

\eeT

See Appendix \ref{Laglogapp} for a proof.

\subsection{Lagrangian Negative Binomial Distribution}

A random variable $X$ is said to have a Lagrangian logarithmic
distribution if it possesses a probability mass function
\[
f(x, \se) = \Pr \{ X = x \} = \f{\ba}{ \al x + \ba} \bi{\al x + \ba}{x} \se^x (1 - \se)^{\ba + \al x - x}, \qqu x = 0, 1, 2, \cd
\]
where $0 < \se < 1, \; \se < \al \se < 1$ and $\ba > 0$.  Let $X_1, \cd, X_n$ be i.i.d. samples of random variable $X$.  By virtue of the LR
method, we have obtained the following result.

\beT

\la{NegBin} \be \la{LagNeg}
 \Pr \{ \ovl{X}_n \leq z \} \leq \li [
\li ( \f{\se}{\vse} \ri )^z \li ( \f{1 - \se}{1 - \vse} \ri )^{\ba +
\al z - z}  \ri ]^n \qu \tx{for} \; 0 \leq z \leq \f{\ba \se}{1 -
\al \se}, \ee where $\vse = \f{z}{\ba + \al z}$. \eeT

See Appendix \ref{NegBinapp} for a proof.

\subsection{Laplace Distribution}

A random variable $X$ is said to have a Lagrangian logarithmic
distribution if it possesses a probability density function
\[
f(x) = \f{1}{2 \ba} \exp \li ( - \f{ | x - \al|  }{\ba}  \ri ), \qqu
- \iy < x < \iy,
\]
where $- \iy < \al < \iy$ and $\ba > 0$.  Let $X_1, \cd, X_n$ be
i.i.d. samples of random variable $X$. By virtue of the LR method,
we have obtained the following results.

\beT

\la{laplace}

\bel &  &  \Pr \{ \ovl{X}_n \geq z \} \leq \li [ \f{z - \al}{\ba}
\exp  \li ( 1 - \f{z - \al}{\ba} \ri ) \ri ]^n \qu \tx{for $z \geq
\al + \ba$},  \la{laplaceineqa}\\
&  & \Pr \{ \ovl{X}_n \leq z \} \leq \li [ \f{\al - z}{\ba} \exp \li
( 1 - \f{\al - z}{\ba} \ri ) \ri ]^n \qu \tx{for $z \leq \al  -
  \ba$}. \la{laplaceineqb}
\eel

\eeT

See Appendix \ref{laplaceapp} for a proof.

\subsection{Logarithmic Distribution}

A random variable $X$ is said to have a logarithmic distribution if
it possesses a probability mass function
\[
f(x) = \f{q^x}{- x \ln p}, \qqu x = 1, 2, \cd
\]
where $p \in (0,1)$ and $q = 1 - p$. Let $X_1, \cd, X_n$ be
i.i.d. samples of random variable $X$.  By virtue of the LR method, we have obtained
the following result.

\beT

\la{ineqlogarithmic}

\be \la{logarithmicsingle}
 \Pr \{ \ovl{X}_n \leq z \} \leq \li [ \f{ \ln (1 - q) }{ \ln (1 - \vse) } \li (
\f{q}{\vse} \ri )^z \ri ]^n \qu \tx{for} \; 0 < z \leq \f{q}{(1 - q)
\ln \f{1}{1 - q}},  \ee where $\vse \in (0, q]$ is the unique number
such that $z = \f{\vse}{(1 - \vse) \ln \f{1}{1 - \vse}}$.

\eeT

See Appendix \ref{ineqlogarithmicapp} for a proof.

\subsection{Lognormal Distribution}

A random variable $X$ is said to have a lognormal distribution if it
possesses a probability density function
\[
f(x) = \f{1}{ x \sq{2 \pi} \si  }  \exp \li [  - \f{1}{2 \si^2} (
\ln x - \mu)^2 \ri ], \qqu x > 0, \qqu -\iy < \mu < \iy, \qqu \si >
0.
 \]
Let $X_1, \cd, X_n$ be i.i.d. samples of random variable $X$.  By
virtue of the LR method, we have obtained the following result.

\beT

\la{lognormal}

\be \la{lognormineq}
 \Pr \{ \ovl{X}_n \leq z \} \leq  \exp \li [ -
\f{n}{2} \li ( \f{\mu - \ln z}{\si} \ri )^2  \ri ] \qu \tx{for $0 <
z \leq e^\mu$}. \ee

\eeT

See Appendix \ref{lognormalapp} for a proof.

\subsection{Nakagami Distribution}

A random variable $X$ is said to have a Nakagami distribution if it
possesses a probability density function
\[
f(x) = \f{2}{\Ga(m)} \f{ x^{2m - 1} }{ \si^{2m} } \exp \li ( - \f{x^2}{\si^2} \ri ),   \qqu x > 0
\]
where $m \geq \f{1}{2}$ and $\si > 0$.  Let $X_1, \cd, X_n$ be
i.i.d. samples of random variable $X$.  By virtue of the LR method,
we have obtained the following results.

\beT

\la{nakagami}

\bel &  &  \Pr \li \{ \ovl{X}_n \leq \f{ \Ga ( \vse + \f{1}{2} ) }{
\Ga (\vse) } \si \ri \} \leq \li \{   \f{ \Ga ( \vse ) }{ \Ga (m) }
\li [ \f{ \Ga ( \vse + \f{1}{2} ) }{ \Ga (\vse) } \ri ]^{2(m-\vse)}
\ri \}^n \qu \tx{for} \; 0 < \vse \leq m, \la{nakagamia} \\
&  & \Pr \{ \ovl{X}_n \geq z \} \leq \li [ \li ( \f{z^2}{m \si^2}
\ri )^{m} \exp \li ( m  - \f{z^2}{\si^2} \ri ) \ri ]^n \qu \tx{for
$z \geq \sq{m} \si$}. \la{nakagamib} \eel

\eeT

See Appendix \ref{nakagamiapp} for a proof.

\subsection{Pareto Distribution}

A random variable $X$ is said to have a Pareto distribution if it
possesses a probability density function
\[
f(x) = \f{\se}{a}  \li ( \f{a}{x} \ri )^{\se + 1}, \qqu x > a > 0,
\qqu \se > 1.
\]
Let $X_1, \cd, X_n$ be i.i.d. samples of random variable $X$.
By virtue of the LR method, we have obtained the following result.

\beT

\la{pareto}

\be \la{paretoineq}
 \Pr \{ \ovl{X}_n \leq \ro \mu \} \leq \li [ e
\se \li ( \f{\se - 1}{\ro \se} \ri )^\se \ln \li ( \f{\ro \se}{\se -
1} \ri ) \ri ]^n \qu \tx{for $1 - \f{1}{\se} < \ro \leq \li ( 1 -
\f{1}{\se} \ri ) \exp \li ( \f{1}{\se} \ri )$}, \ee where $\mu =
\bb{E} [ X ] = \f{ \se a }{ \se - 1 }$.

\eeT

See Appendix \ref{paretoapp} for a proof.

\subsection{Power-Law Distribution}

A random variable $X$ is said to have a power-law distribution if it
possesses a probability density function
\[
f(x) = \f{x^{-\al}}{C(\al)}, \qqu 1 \leq x \leq \ba,
\]
where $\ba > 1, \; \al \in \bb{R}$ and
\[
C(\al) = \bec \f{ 1 - \ba^{1 - \al} }{\al - 1} & \tx{for $\al \neq 1$},\\
\ln \ba &  \tx{for $\al = 1$} \eec
\]
Let $X_1, \cd, X_n$ be i.i.d. samples of random variable $X$.
By virtue of the LR method, we have obtained the following result.

\beT

\la{powerlaw}

 Let $\se \geq \al > 1$ and $z = \f{\se - 1}{\se - 2} \f{ \ba^{\se - 1} - \ba
}{ \ba^{\se - 1} -  1 }$.  Then, \be \la{powineq}
 \Pr \li \{
\ovl{X}_n \leq z  \ri \} \leq \li ( \f{\al - 1}{\se - 1} \f{ 1 -
\ba^{1 -\se} }{ 1 - \ba^{1 - \al} }  z^{\se - \al} \ri )^n. \ee

\eeT

See Appendix \ref{powerlawapp} for a proof.

\subsection{Stirling Distribution}

A random variable is said to have a Stirling distribution if it
possesses a probability mass function
\[
\Pr \{ X = x \} = \f{  m! | s(x, m)| \se^x  }{  x! [- \ln (1 - \se)
]^m }, \qqu 0 < \se < 1, \qu x = m , m+1, \cd,
\]
where $s(x, m)$ is the Stirling number of the first kind, with arguments $x$ and $m$. Let $X_1, \cd, X_n$ be i.i.d. samples of random variable
$X$.  By virtue of the LR method, we have obtained the following result.

\beT

\la{Stirling} \be \la{Stringineq}
 \Pr \li \{ \ovl{X}_n \leq  z  \ri
\} \leq \li [  \f{ \ln (1 - \vse) } {  \ln ( 1 - \se)  } \ri ]^{n m}
\li ( \f{\se}{\vse} \ri )^{n z}  \qu \tx{for} \; z \leq \f{m \se}
{(\se - 1) \ln (1 - \se)}, \ee where $\vse \in (0, \se]$ is the
unique number such that $z = \f{m \vse} {(\vse - 1) \ln (1 -
\vse)}$.

\eeT

See Appendix \ref{Stirlingapp} for a proof.

\subsection{Snedecor's F-Distribution}

If random variable $X$ has a probability density function of the form
\[
f(x) = \f{ \Ga ( \f{n + m}{2}  )  ( \f{m}{n} )^{m \sh 2}  x ^{ (m - 2) \sh 2} }
{ \Ga ( \f{m}{2} ) \Ga (\f{n}{2} )  ( 1 +  \f{m}{n} x )^{(n + m)
\sh 2} }, \qqu \tx{for} \qu 0 < x < \iy,
\]
then the random variable $X$ is said to possess
 an $F$-distribution with $m$ and $n$ degrees of freedom. Making use of the LR method, we have
obtained the following results.

\beT

\la{Snedecor}

\bel &  &  \Pr \{  X \geq z \} \leq z^{m \sh 2}  \li (  \f{ n + m }{
n + m z }  \ri )^{(n + m) \sh 2} \qqu \tx{for} \; z \geq 1
\la{Snedecorineq88a}\\
 &  & \Pr \{  X \leq z \} \leq z^{m \sh 2}
\li ( \f{ n + m }{ n + m z } \ri )^{(n + m) \sh 2} \qqu \tx{for} \;
0 < z \leq 1.  \la{Snedecorineq88b} \eel

\eeT

See Appendix \ref{Snedecorapp} for a proof.

\subsection{Student's t-Distribution}

If random variable $X$ has a probability density function of the form
\[
f(x) = \f{ \Ga ( \f{n + 1}{2}  ) } { \sq{n \pi} \Ga (\f{n}{2} )  ( 1
+ \f{x^2}{n}  )^{(n + 1) \sh 2} }, \qqu \tx{for} \qu - \iy < x <
\iy,
\]
then the random variable $X$ is said to possess
 a Student's $t$-distribution with $n$ degrees of freedom. By virtue of the LR method, we have
obtained the following results.

\beT

\la{CDFineq}

\bel &  &  \Pr \{  | X | \geq z \} \leq z \li (  \f{ n + 1}{ n + z^2
} \ri )^{(n + 1) \sh 2} \qqu \tx{for} \; z \geq 1,  \la{stuineq88a} \\
&  & \Pr \{  | X | \leq z \} \leq z \li (  \f{ n + 1}{ n + z^2 } \ri
)^{(n + 1) \sh 2} \qqu \tx{for} \; 0 < z \leq 1. \la{stuineq88b}
\eel

\eeT

See Appendix \ref{CDFineqapp} for a proof.

\subsection{Truncated Exponential Distribution}

A random variable $X$ is said to have a truncated exponential
distribution if it possesses a probability density function
\[
f(x) = \f{ \se e^{\se x} }{ e^\se - 1 }, \qqu \se \neq 0, \qqu 0 < x
< 1.
\]
Let $X_1, \cd, X_n$ be i.i.d. samples of random variable $X$.  By
virtue of the LR method, we have obtained the following results.

\beT

\la{truncateexp}

\be \la{expch}
 \Pr \{ \ovl{X}_n \leq  z \} \leq \li [  \f{ \se }{ \vse } \f{ e^\vse - 1
}{ e^\se - 1 } e^{ (\se - \vse) z } \ri ]^n  \qu \tx{for}  \; 0 < z
\leq 1 + \f{1}{e^\se - 1} - \f{1}{\se} \; \tx{and} \; z \neq
\f{1}{2}, \ee where $\vse \in (-\iy , \se], \; \vse \neq 0$
satisfies equation $z = 1 + \f{1}{e^\vse - 1} - \f{1}{\vse}$.
Moreover, \be \la{expch88b}
 \Pr \li \{ \ovl{X}_n \leq  \f{1}{2} \ri
\} \leq \li (  \f{ \se e^{\se \sh 2}  }{ e^\se - 1 }  \ri )^n \qu
\tx{for} \; \se > 0. \ee

\eeT

See Appendix \ref{truncateexpapp} for a proof.

\subsection{Uniform Distribution}

Let $X$ be a random variable uniformly distributed over interval $[0, 1]$.
Let $X_1, \cd, X_n$ be i.i.d. samples of the random variable $X$.
By virtue of the LR method, we have obtained the following results.

\beT

\la{UniformCDF}

 \be \la{uniformineqa}
 \Pr \{ \ovl{X}_n \geq z \} \leq \li
( \f{ e^\vse - 1 }{ \vse e^{ \vse z } }   \ri )^n \leq \exp \li ( -
6 n \li ( z - \f{1}{2} \ri )^2 \ri ) \qqu \tx{for} \; 1 > z >
\f{1}{2}, \ee where $\vse$ is a positive number such that $z = 1 +
\f{1}{e^\vse - 1} - \f{1}{\vse}$. Similarly, \be
 \la{uniformineqb}
  \Pr \{ \ovl{X}_n \leq z \} \leq \li (  \f{ e^\vse - 1 }
 { \vse e^{ \vse z } }   \ri )^n \leq \exp \li ( -  6 n \li ( z - \f{1}{2}
\ri )^2 \ri )\qqu \tx{for} \; 0 < z < \f{1}{2}, \ee where $\vse$ is
a negative number such that $z = 1 + \f{1}{e^\vse - 1} -
\f{1}{\vse}$.

\eeT

See Appendix \ref{UniformCDFapp} for a proof.

\subsection{Weibull Distribution}

A random variable $X$ is said to have a Weibull distribution if it
possesses a probability density function
\[
f(x) = \al \ba x^{\ba - 1} \exp \li ( - \al x^\ba \ri ), \qqu x > 0,
\qqu \al > 0, \qqu \ba > 0.
\]
Let $X_1, \cd, X_n$ be i.i.d. samples of the random variable $X$.
By virtue of the LR method, we have obtained the following results.

\beT

\la{WeibullCDF}

\bel &  &  \Pr \{ \ovl{X}_n \leq z \} \leq  \li [  \al z^{\ba} \exp (1 - \al z^\ba ) \ri ]^n \qu \tx{for $\al z^\ba \leq 1$ and $\ba < 1$},
\la{Weib88a}\\
 &  &  \Pr \{ \ovl{X}_n \geq z \} \leq  \li [  \al
z^{\ba} \exp (1 - \al z^\ba ) \ri ]^n \qu \tx{for $\al z^\ba \geq 1$ and $\ba > 1$}. \la{Weib88b} \eel

\eeT

See Appendix \ref{WeibullCDFapp} for a proof.

\section{Concentration Inequalities for Multivariate Distributions}

In this section, we shall apply the LR method to derive concentration inequalities
for the joint distributions of multiple random variables.

\subsection{Dirichlet-Compound Multinomial Distribution  }

Random variables $X_1, \cd, X_k$ are said to have a
Dirichlet-compound multinomial distribution if they possess a
probability mass function
\[
f( \bs{x} )  = \bi{n}{ \bs{x} } \f{ \Ga (\sum_{\ell = 0}^k \al_\ell)
}{ \Ga ( n + \sum_{\ell = 0}^k \al_\ell ) }  \prod_{\ell = 0}^k \f{
\Ga (x_\ell + \al_\ell) }{ \Ga (\al_\ell) },  \] where
\[
\bs{x} = [x_0, x_1, \cd, x_k]^\top, \qqu \bi{n}{ \bs{x} } = \f{n!}{ \prod_{\ell = 0}^k x_\ell ! } \]
 and
\[ \sum_{\ell = 0}^k x_\ell = n
\]
with $x_\ell \geq 0$ and $\al_\ell > 0$ for $\ell = 0, 1, \cd, k$.
Based on the LR method, we have obtained the following result.

\beT

\la{DeriCompound}

Assume that $0 < z_\ell \leq \f{n \al_\ell}{\sum_{i = 0}^k \al_i}$
for $\ell = 1, \cd, k$.  Then,  \be \la{DeC88}
 \Pr \{ X_\ell \leq z_\ell, \; \ell = 1, \cd, k \}
\leq   \f{ \Ga (\sum_{\ell = 0}^k \al_\ell) \; \Ga ( n + \sum_{\ell
= 0}^k \se_\ell ) }{ \Ga (\sum_{\ell = 0}^k \se_\ell) \; \Ga ( n +
\sum_{\ell = 0}^k \al_\ell )  }  \prod_{\ell = 1}^k \f{ \Ga (x_\ell
+ \al_\ell) \; \Ga (\se_\ell) }{ \Ga (x_\ell + \se_\ell) \; \Ga
(\al_\ell) }, \ee where $\se_0 = \al_0$ and
\[
\se_\ell = \f{\al_0 z_\ell} {  n - \sum_{i=1}^k z_i }, \qqu \ell =
1, \cd, k.
\]

\eeT

See Appendix \ref{DeriCompoundapp} for a proof.

\subsection{Inverse Matrix Gamma Distribution}

A positive-definite random matrix $\bs{X}$ is said to have an
inverse matrix gamma distribution \cite{Gupta} if it possesses a probability
density function
\[
f(\bs{x}) = \f{| \bs{\Psi} |^{\al}}{ \ba^{p \al}   \Ga_p (\al ) } |
\bs{x} |^{- \al -( p + 1) \sh 2} \exp \li (  - \f{1}{\ba} \tx{tr} (
\bs{\Psi} \bs{x}^{-1} )\ri ),
\]
where $\ba > 0$ is the scale parameter, $\bs{\Psi}$ is a
positive-definite real matrix of size $p \times p$. Here $\bs{x}$ is
a positive-definite matrix of size $p \times p$, and $\Ga_p(.)$ is
the multivariate gamma function.  The inverse matrix gamma
distribution reduces to the Wishart distribution with $\ba = 2, \;
\al = \f{n}{2}$.  Let $\preccurlyeq$ denote the relationship of two
matrices $A$ and $B$ of the same size such that $A \preccurlyeq B$
implies that $B - A$ is positive definite. By virtue of the LR
method, we have obtained the following result.

\beT

\la{Wishart}

\be \la{invGaM}
 \Pr \li \{ \bs{X} \preccurlyeq \ro \bs{\Up} \ri \}
\leq \f{ 1  }{ \ro^{p \al} } \exp \li ( - \f{p}{2} \li ( \f{1}{\ro}
- 1 \ri ) ( 2 \al - p - 1 ) \ri )  \qu \tx{for $0 < \ro < 1$}, \ee
where $\bs{\Up} = \bb{E} [ \bs{X} ] = \f{2}{\ba}  \f{ \bs{\Psi} }{ 2
\al - p - 1 } $ is the expectation of $\bs{X}$.

\eeT

See Appendix \ref{Wishartapp} for a proof.

\subsection{Multivariate Normal Distribution}

A random vector $\bs{X}$ is said to have a multivariate normal
distribution if it possesses a probability density function
\[
f(\bs{x}) = ( 2 \pi)^{- k \sh 2} | \bs{\Si} |^{-1 \sh 2} \exp \li (
- \f{1}{2} (\bs{x} - \bs{\mu})^{\top}
 \bs{\Si}^{-1} (\bs{x} - \bs{\mu})  \ri ),
\]
where $k$ is the dimension of $\bs{X}$,  $\bs{x}$ is a vector of $k$
elements, $\bs{\mu}$ is the expectation of $\bs{X}$, and $\bs{\Si}$
is the covariance matrix of $\bs{X}$.   Let $\bs{X}_1, \cd,
\bs{X}_n$ be i.i.d. samples of $\bs{X}$.  Define
\[
\ovl{\bs{X}}_n = \f{ \sum_{i=1}^n \bs{X}_i }{n}.
\]
Let $\succcurlyeq$ denote the relationship of two vectors $A = [a_1,
\cd, a_k]$ and $B =[b_1, \cd, b_k]$ such that $A \succcurlyeq B$
implies $a_\ell \geq b_\ell, \; \ell = 1, \cd, k$.  By virtue of the
LR method, we have obtained the following result.

\beT

\la{mulnormal}

\be \la{multinormal}
 \Pr \{ \ovl{\bs{X}}_n  \succcurlyeq  \bs{z} \}
\leq \li [ \exp \li ( \bs{\mu}^{\top} \bs{\Si}^{-1} \bs{z}
 - \f{1}{2} [ \bs{z}^{\top} \bs{\Si}^{-1} \bs{z} +   \bs{\mu}^{\top}
  \bs{\Si}^{-1} \bs{\mu} ]  \ri ) \ri ]^n
\ee
 provided that $\bs{\Si}^{-1}  \bs{z} \succcurlyeq \bs{\Si}^{-1}
\bs{\mu}$.

\eeT

See Appendix \ref{mulnormalapp} for a proof.

\subsection{Multivariate Pareto Distribution}

Random variables $X_1, \cd, X_k$ are said to have a multivariate
Pareto distribution if they possess a probability density function
\[
f(x_1, \cd, x_k) = \li ( \prod_{i = 1}^k \f{\al + i - 1}{\ba_i} \ri
)  \li ( 1 - k + \sum_{i=1}^k \f{x_i}{\ba_i} \ri )^{- (\al + k) },
\qqu x_i
> \ba_i
> 0, \qu \al > 0.
\]
Let $\bs{\fra{X}} = [X_1, \cd, X_k]^\top$. Let $\bs{z} = [z_1, \cd,
z_k]^\top$.  Let $\bs{\fra{X}}_1, \cd, \bs{\fra{X}}_n$ be i.i.d.
samples of random vector $\bs{\fra{X}}$. Define \[ \ovl{
\bs{\fra{X}} }_n = \f{ \sum_{i=1}^n \bs{\fra{X}}_i }{n}.
\]
Let the notation ``$\preceq$'' denote the relationship of two
vectors $A = [a_1, \cd, a_k]^\top$ and $B = [b_1, \cd, b_k]^\top$
such that $A \preceq B$ means $a_\ell \leq b_, \; \ell = 1, \cd, k$.

By virtue of the LR method, we have the following results.

\beT

\la{multiPareto} Let $z_\ell > \ba_\ell, \; \ell = 1, \cd, k$.  The
following statements hold true.

(I):  The inequality \be \la{mulparetoineq}
 \Pr \{ \ovl{
\bs{\fra{X}} }_n \preceq \bs{z} \} \leq  \li [  \li ( \prod_{i =
1}^k \f{\al + i - 1}{\se + i - 1} \ri )
 \li ( 1 - k + \sum_{i=1}^k \f{z_i}{\ba_i} \ri )^{\se - \al} \ri
 ]^n
\ee holds for any $\se > \al$.

(II): The inequality (\ref{mulparetoineq}) holds for $\se$ such that
\be \la{bestp88}
 \sum_{\ell = 0}^{k-1} \f{1}{\se + \ell}  = \ln \li
( 1 - k + \sum_{i=1}^k \f{z_i}{\ba_i} \ri ) \ee provided that \be
\la{best9988}
 \sum_{\ell = 0}^{k-1} \f{1}{\al + \ell}  > \ln \li ( 1 -
k + \sum_{i=1}^k \f{z_i}{\ba_i} \ri ). \ee

(III): The inequality (\ref{mulparetoineq}) holds for \be
\la{mean889}
 \se = 1 + \f{1}{ \li ( \f{1}{k}  \sum_{i=1}^k
\f{z_i}{\ba_i} \ri ) - 1 } \ee provided that $\al > 1$ and $\f{1}{k}
\sum_{i=1}^k \f{z_i}{\ba_i} < \f{\al}{\al - 1}$.

\eeT

See Appendix \ref{multiParetoapp} for a proof.

\section{Conclusion}

We have investigated the concentration phenomenon of random variables based on the likelihood ratio method.  A wide variety of concentration
inequalities for various distributions are developed without using moment generating functions.  The new inequalities are generally simple,
insightful and fairy tight.

\appendix

\section{Proofs of Univariate Inequalities}

\subsection{Proof of Theorem \ref{ChenChBE} } \la{ChenChBEapp}

Let $f(.)$ denote the pmf or pdf of random variable $X$.  Let
$\varTheta$ be the set of non-negative real number such that the
moment generating function $\phi (.)$ of $X$ exists.  Define
\[
g(x, \vse) = \f{ f(x) e^{\vse x} }{ \phi(\vse)}, \qqu \vse \in
\varTheta.
\]
Then, $g(x, \vse)$ is a family of pmf or pdf, which contains $f(x) =
g(x, 0)$. Let $\bs{\mcal{X}} = [X_1, \cd, X_n]$ and $\bs{x} = [x_1,
\cd, x_n]$. The joint pmf or pdf of $\bs{\mcal{X}}$ is
$f_{\bs{\mcal{X}}} (\bs{x}) = \prod_{i=1}^n f(x_i)$, which is
contained in the family
\[
g_{\bs{\mcal{X}}} (\bs{x}, \vse) = \li [ \f{ 1 }{ \phi(\vse)} \ri
]^n \prod_{i=1}^n f(x_i) \exp ( - \vse x_i ), \qqu \fa \vse \in
\varTheta.
\]
It can be checked that
\[
\f{ f_{\bs{\mcal{X}}} (\bs{x}) }{ g_{\bs{\mcal{X}}} (\bs{x}, \vse) }
= \li [ \phi(\vse) \exp \li ( - \vse \ovl{\bs{x}}_n \ri ) \ri ]^n,
\qqu \fa \vse \in \varTheta,
\]
where $\ovl{\bs{x}}_n  = \f{ \sum_{i=1}^n x_i }{n}$. Hence,
\[
\f{ f_{\bs{\mcal{X}}} (\bs{x}) }{ g_{\bs{\mcal{X}}} (\bs{x}, \vse) }
\leq \li [ \phi(\vse) e^{-z}  \ri ]^n, \qqu \fa \vse \in \varTheta
\; \tx{provided that} \; \ovl{\bs{x}}_n \geq z.
\]
This implies that
\[
f_{\bs{\mcal{X}}} (\bs{\mcal{X}}) \; \bb{I}_{ \{  \ovl{X}_n \geq z \} } \leq \Lm (\vse) \; g_{\bs{\mcal{X}}} (\bs{\mcal{X}}, \vse),
\]
where $\Lm (\vse) = \li [ \phi(\vse) e^{-z}  \ri ]^n$. By
differentiation, it can be shown that the infimum of $\Lm (\vse)$
with respect to $\vse \in \varTheta$ is attained at $\tau \in
\varTheta$ such that $\f{ \phi^\prime (\tau) }{\phi (\tau) } = z$.
It follows from (\ref{maintight}) of Theorem \ref{ThM888} that \be
\la{vip96638} \Pr \{ \ovl{X}_n \geq z \} \leq \Lm (\tau) \;
\bb{P}_\tau \{ \ovl{X}_n \geq z \} \leq \Lm (\tau) = \li [
\phi(\tau) e^{-z \tau}  \ri ]^n. \ee

Now we evaluate $\bb{P}_\tau \{ \ovl{X}_n \geq z \}$.  Let $\bb{E}_\tau [. ]$ denote
the expectation of a function of random variable $X$ having
pmf or pdf $g(x, \tau)$. Note that
\[
\bb{E}_\tau [ X ] = \int \f{ x f(x) e^{\tau x} }{\phi(\tau)} dx =
\f{1}{ \phi(\tau) } \int x f(x) e^{\tau x} dx =
\f{\phi^\prime(\tau)}{ \phi(\tau) } = z.
\]
Similarly,
\[
\bb{E}_\tau [ X^2 ] = \f{\phi^{\prime \prime} (\tau) }{ \phi (\tau)
}, \qqu  \bb{E}_\tau [ X^3 ] = \f{\phi^{\prime \prime \prime} (\tau)
}{ \phi (\tau) }, \qqu  \bb{E}_\tau [ X^4 ] = \f{\phi^{\prime \prime
\prime \prime} (\tau) }{ \phi (\tau) }
\]
So,
\[
\bb{E}_\tau [ |X - z|^2 ] = \bb{E}_\tau [ X^2 ] - z^2 = \f{\phi^{\prime \prime} (\tau) }{ \phi (\tau) } - z^2.
\]
Note that $(X - z)^4 = X^4 - 4 z X^3 + 6 z^2 X^2 - 4 z^3 X + z^4$. Hence, \bee  \bb{E}_\tau [ (X - z)^4 ]  =  \f{1}{ \phi (\tau) }  \li [
\phi^{\prime \prime \prime \prime} (\tau) - 4 z  \phi^{ \prime \prime \prime } (\tau) + 6 z^2 \phi^{ \prime \prime } (\tau)  - 3 z^4 \phi (\tau)
\ri ].  \eee From Berry-Essen's inequality \cite{Berry, Essen}, we have \bee  \bb{P}_\tau \{ \ovl{X}_n \geq z \} & \leq &  \li \{  \f{1}{2} +
\f{ C_{BE} }{\sq{n}} \li [ \f{ \phi^{\prime \prime \prime \prime} (\tau) - 4 z  \phi^{ \prime \prime \prime } (\tau) + 6 z^2 \phi^{ \prime
\prime } (\tau) - 3 z^4 \phi (\tau) } { ( \phi^{\prime \prime} (\tau) -  z^2 \phi(\tau)  )^2 \sh \phi (\tau) } \ri
]^{\f{3}{4}} \ri \}\\
& = & \li \{  \f{1}{2} + \f{ C_{BE} }{\sq{n}} \li ( \f{ \phi (\tau)
[ \phi^{\prime \prime \prime \prime} (\tau) - 4 z  \phi^{ \prime
\prime \prime } (\tau) ] + 3 [ \phi^{ \prime \prime } (\tau) ]^2 } {
[ \phi^{\prime \prime} (\tau) - z^2 \phi(\tau)  ]^2  } - 3 \ri
)^{\f{3}{4}} \ri \}.  \eee Making use of the above inequalities and
(\ref{vip96638}) completes the proof of the theorem.

\subsection{Proof of Theorem \ref{betathm} } \la{betathmapp}

Let $\bs{\mcal{X}} = [ X_1, \cd, X_n ]$ and $\bs{x} = [x_1, \cd,
x_n]$.  The joint
 probability density function of $\bs{\mcal{X}}$ is
 \[
f_{\bs{\mcal{X}}} (\bs{x}) = \frac{1}{ [\mcal{B} (\al, \ba)]^n } \li
( \prod_{i=1}^n x_i \ri )^{\al - 1} \li [ \prod_{i=1}^n (1 - x_i)
\ri ]^{\ba - 1}.
 \]
To apply the LR method, we construct a family of probability density functions
\[
g_{\bs{\mcal{X}}} (\bs{x}, \vse) = \frac{1}{ [\mcal{B} (\vse,
\ba)]^n } \li ( \prod_{i=1}^n x_i \ri )^{\vse - 1} \li [
\prod_{i=1}^n (1 - x_i) \ri ]^{\ba - 1}
 \]
for $\vse \in (0, \al]$.  It can be checked that
\[
\f{ f_{\bs{\mcal{X}}} (\bs{x}) }{ g_{\bs{\mcal{X}}} (\bs{x}, \vse) }
= \li [  \f{ \mcal{B} (\vse, \ba) } {  \mcal{B} (\al, \ba) }  \ri
]^n \li ( \prod_{i=1}^n x_i \ri )^{\al - \vse}.
\]
Since the geometric mean is no greater than the arithmetic mean, we
have
\[
\prod_{i=1}^n x_i  \leq \li (  \ovl{\bs{x}}_n  \ri )^n,
\]
where $\ovl{\bs{x}}_n  = \f{ \sum_{i=1}^n x_i }{n}$.  Hence,
\[
\f{ f_{\bs{\mcal{X}}} (\bs{x}) }{ g_{\bs{\mcal{X}}} (\bs{x}, \vse) }
\leq \li [  \f{ \mcal{B} (\vse, \ba) } {  \mcal{B} (\al, \ba) } \li
(  \ovl{\bs{x}}_n  \ri )^{\al - \vse} \ri ]^n
\]
and it follows that
\[
\f{ f_{\bs{\mcal{X}}} (\bs{x}) }{ g_{\bs{\mcal{X}}} (\bs{x}, \vse) } \leq \li [  \f{ \mcal{B} (\vse, \ba) } {  \mcal{B} (\al, \ba) } z^{\al -
\vse} \ri ]^n \qu \fa \vse \in (0, \al] \; \tx{provided that} \; \ovl{\bs{x}}_n \leq z.
\]
Consequently,
\[
f_{\bs{\mcal{X}}} ( \bs{\mcal{X}}) \; \bb{I}_{ \{ \ovl{X}_n \leq z \} } \leq \Lm (\vse) \; g_{\bs{\mcal{X}}} ( \bs{\mcal{X}}, \vse) \qu \fa \vse
\in (0, \al],
\]
where
\[
\Lm (\vse) = \li [  \f{ \mcal{B} (\vse, \ba) } {  \mcal{B} (\al,
\ba) } z^{\al - \vse} \ri ]^n.
\]
It follows from Theorem \ref{ThM888} that \be \la{recall} \Pr \{ \ovl{X}_n \leq z \}  \leq  \li [ \f{ 1 } {  \mcal{B} (\al, \ba) } \; \inf_{\vse
\in (0, \al] }  \mcal{B} (\vse, \ba)   z^{\al - \vse} \ri ]^n \qu \tx{for $0 < z < 1$}. \ee
  As a consequence of $0 < z \leq \mu$ and the
definition of $\wh{\al}$, we have that $0 < \wh{\al} \leq \al$.  Hence,
\[
\inf_{\vse \in (0, \al] }  \mcal{B} (\vse, \ba)   z^{\al - \vse}
\leq  \mcal{B} (\wh{\al}, \ba)   z^{\al - \wh{\al}},
\]
which leads to (\ref{beta888a}).

To show (\ref{beta888b}),  we construct a family of probability density functions
\[
g_{\bs{\mcal{X}}} (\bs{x}, \vse) = \frac{1}{ [\mcal{B} (\al,
\vse)]^n } \li ( \prod_{i=1}^n x_i \ri )^{\al - 1} \li [
\prod_{i=1}^n (1 - x_i) \ri ]^{\vse - 1}
 \]
for $\vse \in  (0, \ba]$.  It can be checked that
\[
\f{ f_{\bs{\mcal{X}}} (\bs{x}) }{ g_{\bs{\mcal{X}}} (\bs{x}, \vse) } = \li [  \f{ \mcal{B} (\al, \vse) } {  \mcal{B} (\al, \ba) }  \ri ]^n \li [
\prod_{i=1}^n (1 - x_i) \ri ]^{\ba - \vse} \leq \li [  \f{ \mcal{B} (\al, \vse) } {  \mcal{B} (\al, \ba) } (1 - z)^{\ba - \vse} \ri ]^n \qu \fa
\vse \in  (0, \ba] \; \tx{provided that} \; \ovl{\bs{x}}_n \geq z.
\]
Hence,
\[
f_{\bs{\mcal{X}}} ( \bs{\mcal{X}}) \; \bb{I}_{ \{ \ovl{X}_n \geq z \} } \leq \Lm (\vse) \; g_{\bs{\mcal{X}}} ( \bs{\mcal{X}}, \vse) \qu \fa \vse
\in (0, \ba],
\]
where
\[
\Lm (\vse) = \li [  \f{ \mcal{B} (\al, \vse) } {  \mcal{B} (\al,
\ba) } (1 - z)^{\ba - \vse} \ri ]^n.
\]
It follows from Theorem \ref{ThM888} that \bee \Pr \{ \ovl{X}_n \geq z \} \leq \li [ \f{ 1 } { \mcal{B} (\al, \ba) } \; \inf_{\vse \in (0, \ba]
} \mcal{B} (\al, \vse) (1 - z)^{\ba - \vse} \ri ]^n \qu \tx{for $0 < z < 1$}. \eee As a consequence of $\mu \leq z < 1$ and the definition of
$\wh{\ba}$, we have that $0 < \wh{\ba} \leq \ba$. Hence,
\[
\inf_{\vse \in (0, \ba] }  \mcal{B} (\al, \vse)   (1 - z)^{\ba -
\vse} \leq \mcal{B} (\al, \wh{\ba})   (1 - z)^{\ba - \wh{\ba}},
\]
which leads to (\ref{beta888b}).

Finally, we need to show (\ref{beta888c}). Since $\ba = 1$,  using $\Ga(z + 1) = z \Ga(z)$, we obtain from (\ref{recall}) the following
inequality \be \la{reuse88a} \Pr \li \{ \ovl{X}_n \leq z \ri  \} \leq \li [ \inf_{\vse \in (0, \al] } \f{ \al z^{\al - \vse} } { \vse } \ri ]^n.
\ee Consider function $w (\vse) = \ln \al - \ln \vse + (\al - \vse) \ln z$. Note that the first and second derivatives are $w^\prime (\vse) = -
\f{1}{\vse} - \ln z$ and $w^{\prime \prime} (\vse) = \f{1}{\vse^2}$, respectively.   By the assumption that $0 < z < \exp ( - \f{1}{\al} )$, the
infimum is attained at $\vse = \f{1}{\ln \f{1}{z} } \in (0, \al)$.  Hence, \be \la{reuse88b} \inf_{\vse \in (0, \al] } w (\vse) = \ln \al - \ln
\f{1}{\ln \f{1}{z} } + \li ( \al - \f{1}{\ln \f{1}{z} } \ri ) \ln z = 1 + \ln \li ( \al z^\al \ln \f{1}{z} \ri ). \ee Combining (\ref{reuse88a})
and (\ref{reuse88b}) yields
\[
\Pr \li \{ \ovl{X}_n \leq z \ri  \} \leq \li [ \exp \li ( 1 + \ln
\li ( \al z^\al \ln \f{1}{z} \ri ) \ri ) \ri ]^n = \li ( e \al z^\al
\ln \f{1}{z} \ri )^n \qu \tx{for $0 < z < \exp \li ( - \f{1}{\al}
\ri )$}.
\]
This proves (\ref{beta888c}).  The proof of the theorem is thus
completed.

\subsection{Proof of Theorem \ref{ineqthm4} } \la{ineqthm4app}

To apply the LR method, we construct a family of probability mass functions
\[
g (x, \vse) = \bi{n + x - 1}{x} \f{ \Ga ( \al + n) \Ga(\vse + x) \Ga
(\al + \vse)  }
 {  \Ga (\al + \vse + n + x) \Ga(\al) \Ga(\vse) }, \qqu x =
0, 1, 2, \cd
 \]
for $\vse \in (0, \ba]$.  Define
\[
L (x, \vse) = \f{ f(x) }{ g (x, \vse)  }, \qqu x = 0, 1, 2, \cd.
\]
Then,
\[
L (x, \vse) = \f{ \Ga(\vse)  }{ \Ga(\ba) } \f{ \Ga(\ba + x) } {
\Ga(\vse + x)  } \f{ \Ga (\al + \vse + n + x) } { \Ga (\al + \ba + n
+ x) }, \qqu x = 0, 1, 2, \cd.
\]
It can be checked that
\[
\f{ L (x+1, \vse) }{ L (x, \vse)} = \f{ \ba + x } { \vse + x } \f{
\al + \vse + n + x } { \al + \ba + n + x  } \geq 1, \qqu x = 0, 1,
2, \cd
\]
for $\vse \in (0, \ba]$.  This implies that for any non-negative
integer $z$,
\[
L (x, \vse)  \leq  L (z, \vse), \qqu \fa \vse \in (0, \ba]
\]
for any non-negative integer $x$ no greater than $z$.  Hence,
\[
\f{ f(x) }{ g (x, \vse)  } \leq \Lm (\vse), \qqu \fa \vse \in (0,
\ba]
\]
for any non-negative integer $x$ no greater than $z$, where
\[
\Lm (\vse) = \f{ \Ga(\vse)  }{ \Ga(\ba) } \f{ \Ga(\ba + z) } {
\Ga(\vse + z)  } \f{ \Ga (\al + \vse + n + z) } { \Ga (\al + \ba + n
+ z) }.
\]
Consequently,
\[
f ( X ) \; \bb{I}_{ \{ X \leq z \} } \leq \Lm (\vse) \; g ( X, \vse) \qu \fa \vse \in (0, \ba].
\]
By virtue of Theorem \ref{ThM888}, we have
\[
\Pr \{ X \leq z \} \leq \inf_{\vse \in (0, \ba]} \Lm (\vse).
\]
Since $z \leq \bb{E}[X] = \f{ n \ba }{\al - 1 }$, we have
\[
0 \leq \f{ \al z - z }{n}  \leq \ba
\]
and thus
\[
\Pr \{ X \leq z \} \leq  \Lm \li ( \f{ \al z - z }{n} \ri ) = \f{
\Ga( \f{ \al z - z }{n} )  }{ \Ga(\ba) } \f{ \Ga(\ba + z) } {
\Ga(\f{ \al z - z }{n} + z)  } \f{ \Ga (\al + \f{ \al z - z }{n} + n
+ z) } { \Ga (\al + \ba + n + z) }.
\]
This completes the proof of the theorem.

\subsection{Proof of Theorem \ref{betaprime}} \la{betaprimeapp}

Let $\bs{\mcal{X}} = [ X_1, \cd, X_n ]$ and $\bs{x} = [x_1, \cd,
x_n]$.  The joint
 probability density function of $\bs{\mcal{X}}$ is
 \[
f_{\bs{\mcal{X}}} (\bs{x}) =  \f{\prod_{i=1}^n [ x_i^{\al - 1} (1 +
x_i)^{-\al - \ba} ]}{ [ \mcal{B}(\al, \ba) ]^n }.
 \]
To apply the LR method to show (\ref{betaprimea}), we construct a family of probability density functions
\[
g_{\bs{\mcal{X}}} (\bs{x}, \vse) =  \f{\prod_{i=1}^n [ x_i^{\vse -
1} (1 + x_i)^{- \vse - \ba} ]}{ [ \mcal{B}(\vse, \ba) ]^n }, \qqu
\vse \in (0, \al].
\]
It can be checked that
\[
\f{ f_{\bs{\mcal{X}}} (\bs{x}) }{g_{\bs{\mcal{X}}} (\bs{x}, \vse)} =
\li [ \f{ \mcal{B}(\vse, \ba) } { \mcal{B}(\al, \ba) } \ri ]^n
\prod_{i=1}^n  \li ( \f{x_i}{1 + x_i} \ri )^{\al - \vse}, \qqu \vse
\in (0, \al].
\]
By differentiation, it can be shown that $\ln \f{x}{1 + x}$ is a
concave function of $x > 0$. As a consequence of this fact, we have
\[
\prod_{i=1}^n  \li ( \f{x_i}{1 + x_i} \ri )^{\al - \vse} \leq \li (
\f{\ovl{\bs{x}}_n}{1 + \ovl{\bs{x}}_n} \ri )^{n(\al - \vse)} \qqu
\fa \vse \in (0, \al],
\]
where $\ovl{\bs{x}}_n = \f{ \sum_{i=1}^n x_i  }{n}$. Since $\f{x}{1
+ x}$ is an increasing function of $x > 0$, it follows that
\[
\prod_{i=1}^n  \li ( \f{x_i}{1 + x_i} \ri )^{\al - \vse} \leq \li (
\f{z}{1 + z} \ri )^{n(\al - \vse)} \qqu \fa \vse \in (0, \al] \qu
\tx{provided that} \; 0 \leq \ovl{\bs{x}}_n \leq z.
\]
Therefore, we have established that
\[
f_{\bs{\mcal{X}}} ( \bs{\mcal{X}}) \; \bb{I}_{ \{ \ovl{X}_n \leq z
\} } \leq \Lm (\vse) \; g_{\bs{\mcal{X}}} ( \bs{\mcal{X}}, \vse)
\qqu \fa \vse \in (0, \al],
\]
where
\[
\Lm (\vse) = \li [ \f{ \mcal{B}(\vse, \ba) } { \mcal{B}(\al, \ba) }
\li ( \f{z}{1 + z} \ri )^{\al - \vse}  \ri ]^n.
\]
Invoking Theorem \ref{ThM888}, we have
\[
\Pr \li \{ \ovl{X}_n \leq z \ri \} \leq \inf_{ \vse \in (0, \al] }
\Lm (\vse).
\]
As a consequence of $\ba > 1$ and $0 < z \leq \f{\al}{\ba - 1}$, we
have $0 < z (\ba - 1) \leq \al$. Hence,
\[
\Pr \li \{ \ovl{X}_n \leq z \ri \} \leq  \Lm ( \ba z - z )  = \li [
\f{ \mcal{B}(\ba z - z, \ba) } { \mcal{B}(\al, \ba) } \li ( \f{z}{1
+ z} \ri )^{\al + z - \ba z} \ri ]^n.
\]
This proves (\ref{betaprimea}).

To apply the LR method to show (\ref{betaprimeb}), we construct a family of probability density functions
\[
g_{\bs{\mcal{X}}} (\bs{x}, \vse) =  \f{\prod_{i=1}^n [ x_i^{\al - 1}
(1 + x_i)^{- \al - \vse} ]}{ [ \mcal{B}(\al, \vse) ]^n }, \qqu \vse
\in [\ba, \iy).
\]
It can be seen that
\[
\f{ f_{\bs{\mcal{X}}} (\bs{x}) }{g_{\bs{\mcal{X}}} (\bs{x}, \vse)} =
\li [ \f{ \mcal{B}(\al, \vse) } { \mcal{B}(\al, \ba) } \ri ]^n
\prod_{i=1}^n  \li ( 1 + x_i \ri )^{\vse - \ba}, \qqu \vse \in [\ba,
\iy).
\]
By differentiation, it can be shown that $\ln (1 + x)$ is a concave
function of $x > 0$. As a consequence of this fact, we have
\[
\prod_{i=1}^n  \li ( 1 + x_i \ri )^{\vse - \ba} \leq \li ( 1 +
\ovl{\bs{x}}_n \ri )^{n (\vse - \ba) } \leq (1 + z)^{ n (\vse - \ba)
}, \qqu \vse \in [\ba, \iy)
\]
provided that $0 \leq \ovl{\bs{x}}_n \leq z$. Hence, we have that
$f_{\bs{\mcal{X}}} ( \bs{\mcal{X}}) \; \bb{I}_{ \{ \ovl{X}_n \leq z \} } \leq
\Lm (\vse) \; g_{\bs{\mcal{X}}} ( \bs{\mcal{X}}, \vse)$ holds for any $\vse \in [\ba, \iy)$, where
\[
\Lm (\vse) = \li [ \f{ \mcal{B}(\al, \vse) } { \mcal{B}(\al, \ba) }
\li ( 1 + z \ri )^{ \vse - \ba }  \ri ]^n.
\]
Making use of Theorem \ref{ThM888}, we have $\Pr \li \{ \ovl{X}_n
\leq z \ri \} \leq \inf_{ \vse \geq \ba } \Lm (\vse)$.  As a
consequence of $0 < z \leq \f{\al}{\ba - 1}$, we have $1 +
\f{\al}{z} \geq \ba$.  Hence,
\[
\Pr \li \{ \ovl{X}_n \leq z \ri \} \leq \Lm \li ( 1 + \f{\al}{z} \ri
) = \li [  \f{ \mcal{B}(\al, 1 + \f{\al}{z}) }{ \mcal{B}(\al, \ba) }
(1 + z)^{1 + \f{\al}{z} - \ba} \ri ]^n \qu \tx{for $0 < z \leq
\f{\al}{\ba - 1}$}.
\]
This proves (\ref{betaprimeb}). The proof of the theorem is thus
completed.

\subsection{Proof of Theorem \ref{borelcdf}} \la{borelcdfapp}

Let $\bs{\mcal{X}} = [ X_1, \cd, X_n ]$ and $\bs{x} = [x_1, \cd,
x_n]$.  The joint
 probability mass function of $\bs{\mcal{X}}$ is
 \[
f_{\bs{\mcal{X}}} (\bs{x}) =  \prod_{i=1}^n \f{(\se x_i)^{x_i - 1}
e^{- \se x_i}}{x_i!}.
 \]
To apply the LR method to show (\ref{borelineq}), we construct a family of probability mass functions
 \[
g_{\bs{\mcal{X}}} (\bs{x}, \vse) =  \prod_{i=1}^n \f{(\vse x_i)^{x_i
- 1} e^{- \vse x_i}}{x_i!}, \qqu \vse \in (0, \se].
 \]
It can be seen that
\[
\f{ f_{\bs{\mcal{X}}} (\bs{x}) }{  g_{\bs{\mcal{X}}} (\bs{x}, \vse)
} = \li [ \li ( \f{\se}{\vse} \ri )^{\ovl{\bs{x}}_n  - 1} \exp \li (
(\vse - \se)  \ovl{\bs{x}}_n \ri ) \ri ]^n = \li [ \li (
\f{\vse}{\se} \ri ) \exp \li ( (\ln \se - \se - \ln \vse + \vse ) \;
\ovl{\bs{x}}_n \ri ) \ri ]^n ,
\]
where $\ovl{\bs{x}}_n = \f{\sum_{i=1}^n x_i}{n}$. Noting that $\ln x
- x$ is increasing with respect to $x \in (0, 1)$, we have that
\[
\ln \se - \se - \ln \vse + \vse \geq 0
\]
as a consequence of $0 < \vse \leq \se$.  It follows that
\[
\li [ \li ( \f{\vse}{\se} \ri ) \exp \li ( (\ln \se - \se - \ln \vse + \vse ) \; \ovl{\bs{x}}_n \ri ) \ri ]^n \leq \li [ \li ( \f{\vse}{\se} \ri
) \exp \li ( (\ln \se - \se - \ln \vse + \vse ) \; z \ri ) \ri ]^n \qu \fa \vse \in (0, \se]
\]
provided that $\ovl{\bs{x}}_n \leq z$.  Hence,
\[
\f{ f_{\bs{\mcal{X}}} (\bs{x}) }{  g_{\bs{\mcal{X}}} (\bs{x}, \vse) } \leq \Lm (\vse) \qu \fa \vse \in (0, \se] \; \tx{provided that} \;
\ovl{\bs{x}}_n \leq z,
\]
where
\[
\Lm (\vse) = \li [ \li ( \f{\vse}{\se} \ri ) \exp \li ( (\ln \se -
\se - \ln \vse + \vse ) \; z \ri ) \ri ]^n.
\]
Hence, we have that $f_{\bs{\mcal{X}}} ( \bs{\mcal{X}}) \; \bb{I}_{
\{ \ovl{X}_n \leq z \} } \leq \Lm (\vse) \; g_{\bs{\mcal{X}}} (
\bs{\mcal{X}}, \vse)$ holds for any $\vse \in (0, \se]$. By virtue
of Theorem \ref{ThM888}, we have $\Pr \li \{ \ovl{X}_n \leq z \ri \}
\leq \inf_{ \vse \in (0, \se] } \Lm (\vse)$.  By differentiation, it
can be shown that the infimum of $\Lm (\vse)$ with respective to
$\vse \in (0, \se]$ is attained at $\vse = 1 - \f{1}{z}$. Therefore,
\[
\Pr \{ \ovl{X}_n \leq z \} \leq \Lm \li ( 1 - \f{1}{z} \ri ) = \li [
\li ( \f{ e \se z }{ 1 - z } \ri )^{z - 1} e^{- \se z} \ri ]^n \qu
\tx{for $1 < z < \f{1}{1 - \se}$}.
\]
This completes the proof of the theorem.

\subsection{Proof of Theorem \ref{consulcdf}} \la{consulcdfapp}

Let $\bs{\mcal{X}} = [ X_1, \cd, X_n ]$ and $\bs{x} = [x_1, \cd,
x_n]$.  The joint
 probability mass function of $\bs{\mcal{X}}$ is
 \[
f_{\bs{\mcal{X}}} (\bs{x}) =  \prod_{i=1}^n \f{1}{x_i} \bi{m
x_i}{x_i -1} \li ( \f{\se}{1 - \se} \ri )^{x_i - 1} (1 - \se)^{m
x_i}.
 \]
To apply the LR method to show (\ref{consineq}), we construct a family of probability mass functions
 \[
g_{\bs{\mcal{X}}} (\bs{x}, \vse) =  \prod_{i=1}^n \f{1}{x_i} \bi{m
x_i}{x_i -1} \li ( \f{\vse}{1 - \vse} \ri )^{x_i - 1} (1 - \vse)^{m
x_i}, \qqu \vse \in (0, \se].
 \]
It can be verified that
\[
\f{ f_{\bs{\mcal{X}}} (\bs{x})  }{ g_{\bs{\mcal{X}}} (\bs{x}, \vse)
} =  \li \{ \li [ \f{\se (1 - \vse)}{\vse (1 - \se) } \ri
]^{\ovl{\bs{x}}_n - 1} \li ( \f{ 1 - \se  } { 1 - \vse } \ri )^{m
\ovl{\bs{x}}_n} \ri \}^n,
\]
where $\ovl{\bs{x}}_n = \f{\sum_{i=1}^n x_i}{n}$. Define function
\[
h (x) = \ln  \f{x}{1 - x} + m \ln (1 - x)
\]
for $x \in (0, 1)$.  Then,
\[
\f{ f_{\bs{\mcal{X}}} (\bs{x})  }{ g_{\bs{\mcal{X}}} (\bs{x}, \vse)
} =  \li \{ \f{\vse (1 - \se) } {\se (1 - \vse)}   \exp \li (
\ovl{\bs{x}}_n [ h(\se) - h (\vse) ] \ri ) \ri \}^n.
\]
Note that the first derivative of $h(x)$ is $h ^\prime (x) = \f{1}{1
- x} \li ( \f{1}{x} - m \ri )$, which is positive for $x \in (0,
1)$.  Hence, $h(\se) - h (\vse) \geq 0$ for $\vse \in (0, \se]$. It
follows that
\[
\f{ f_{\bs{\mcal{X}}} (\bs{x})  }{ g_{\bs{\mcal{X}}} (\bs{x}, \vse)
} \leq  \Lm (\vse) \qu \fa \vse \in (0, \se] \; \tx{provided that}
\; 1 \leq \ovl{\bs{x}}_n \leq z,
\]
where \[ \Lm (\vse) = \li \{ \f{\vse (1 - \se) } {\se (1 - \vse)}
\exp \li ( z [ h(\se) - h (\vse) ] \ri ) \ri \}^n.
\]
This implies that $f_{\bs{\mcal{X}}} ( \bs{\mcal{X}}) \; \bb{I}_{ \{
\ovl{X}_n \leq z \} } \leq \Lm (\vse) \; g_{\bs{\mcal{X}}} (
\bs{\mcal{X}}, \vse)$ holds for any $\vse \in (0, \se]$.  By virtue
of Theorem \ref{ThM888}, we have $\Pr \li \{ \ovl{X}_n \leq z \ri \}
\leq \inf_{ \vse \in (0, \se] } \Lm (\vse)$.  By differentiation, it
can be shown that the infimum of $\Lm (\vse)$ with respective to
$\vse \in (0, \se]$ is attained at $\vse = \f{z - 1}{m z}$.  So,
\[
\Pr \{ \ovl{X}_n \leq z \} \leq  \Lm \li ( \f{z - 1}{m z} \ri ) =
\li [ \f{ \li ( \f{\se}{1 - \se} \ri )^{z - 1} (1 - \se)^{m z}  }
 {  \li (  \f{z - 1}{1 - z + m z} \ri )^{z - 1}
 (1 - \f{z - 1}{m z})^{m z}  } \ri ]^n \qu \tx{for $1 \leq z < \f{1}{1 - m \se}$}.
\]
This completes the proof of the theorem.

\subsection{Proof of Theorem \ref{ineqgeeta}} \la{ineqgeetaapp}

Let $\bs{\mcal{X}} = [ X_1, \cd, X_n ]$ and $\bs{x} = [x_1, \cd,
x_n]$.  The joint
 probability mass function of $\bs{\mcal{X}}$ is
 \[
f_{\bs{\mcal{X}}} (\bs{x}) =  \prod_{i=1}^n \f{1}{\ba x_i - 1}
\bi{\ba x_i - 1}{x_i} \se^{x_i - 1} (1 - \se)^{(\ba - 1) x_i }.
 \]
To apply the LR method to show (\ref{geetaineq}), we construct a family of probability mass functions
 \[
g_{\bs{\mcal{X}}} (\bs{x}, \vse) =  \prod_{i=1}^n \f{1}{\ba x_i - 1}
\bi{\ba x_i - 1}{x_i} \vse^{x_i - 1} (1 - \vse)^{(\ba - 1) x_i },
\qqu \vse \in (0, \se].
 \]
It can be verified that
\[
\f{ f_{\bs{\mcal{X}}} (\bs{x})  }{ g_{\bs{\mcal{X}}} (\bs{x}, \vse)
} =  \li \{ \li ( \f{\se }{\vse} \ri )^{\ovl{\bs{x}}_n - 1} \li (
\f{ 1 - \se  } { 1 - \vse } \ri )^{ (\ba - 1) \ovl{\bs{x}}_n} \ri
\}^n,
\]
where $\ovl{\bs{x}}_n = \f{\sum_{i=1}^n x_i}{n}$. Define function
\[
h (x) = \ln x + (\ba - 1) \ln (1 - x)
\]
for $x \in (0, 1)$.  Then,
\[
\f{ f_{\bs{\mcal{X}}} (\bs{x})  }{ g_{\bs{\mcal{X}}} (\bs{x}, \vse)
} =  \li \{ \f{\vse} {\se}   \exp \li ( \ovl{\bs{x}}_n [ h(\se) - h
(\vse) ] \ri ) \ri \}^n.
\]
Note that the first derivative of $h(x)$ is $h^\prime (x) = \f{1 -
\ba x}{x (1 - x)}$, which is positive for $x \in (0, \f{1}{\ba})$.
Hence, $h(\se) - h(\vse) \geq 0$ for $\vse \in (0, \se]$. It follows
that
\[
\f{ f_{\bs{\mcal{X}}} (\bs{x})  }{ g_{\bs{\mcal{X}}} (\bs{x}, \vse) } \leq \Lm (\vse)   \qu \fa \vse \in (0, \se] \; \tx{provided that} \;
\ovl{\bs{x}}_n \leq z,
\]
where
\[
\Lm (\vse) = \li \{ \f{\vse} {\se}   \exp \li ( z [ h(\se) - h
(\vse) ] \ri ) \ri \}^n.
\]
This implies that $f_{\bs{\mcal{X}}} ( \bs{\mcal{X}}) \; \bb{I}_{ \{
\ovl{X}_n \leq z \} } \leq \Lm (\vse) \; g_{\bs{\mcal{X}}} (
\bs{\mcal{X}}, \vse)$ holds for any $\vse \in (0, \se]$.  By virtue
of Theorem \ref{ThM888}, we have $\Pr \li \{ \ovl{X}_n \leq z \ri \}
\leq \inf_{ \vse \in (0, \se] } \Lm (\vse)$.  By differentiation, it
can be shown that the infimum of $\Lm (\vse)$ with respective to
$\vse \in (0, \se]$ is attained at $\vse = \f{z - 1}{\ba z - 1}$.
Therefore,
\[
\Pr \{ \ovl{X}_n \leq z \} \leq \Lm \li (  \f{z - 1}{\ba z - 1} \ri
) = \li [  \f{ \se^{z-1} (1 - \se)^{\ba z - z}  }{ \li ( \f{z -
1}{\ba z - 1} \ri )^{z-1} \li ( 1 - \f{z - 1}{\ba z - 1} \ri )^{\ba
z - z} } \ri ]^n \qu \tx{for} \; 1 \leq z \leq \f{1 - \se}{1 - \ba
\se}.
\]
This completes the proof of the theorem.

\subsection{Proof of Theorem \ref{ineqgumbel}} \la{ineqgumbelapp}

Let $\bs{\mcal{X}} = [ X_1, \cd, X_n ]$ and $\bs{x} = [x_1, \cd,
x_n]$.  The joint
 probability density function of $\bs{\mcal{X}}$ is
 \[
f_{\bs{\mcal{X}}} (\bs{x}) =  \prod_{i=1}^n \f{1}{\ba} \exp \li (
\f{\mu - x_i}{\ba} - \exp \li ( \f{\mu - x_i}{\ba} \ri )  \ri ).
 \]
To apply the LR method to show (\ref{gumbelineq}), we construct a family of probability density functions
 \[
g_{\bs{\mcal{X}}} (\bs{x}, \vse) =  \prod_{i=1}^n \f{1}{\ba} \exp
\li ( \f{\vse - x_i}{\ba} - \exp \li ( \f{\vse - x_i}{\ba} \ri )
\ri ), \qqu \vse \in (- \iy, \mu].
 \]
Note that \bee \f{ f_{\bs{\mcal{X}}} (\bs{x})  }{ g_{\bs{\mcal{X}}}
(\bs{x}, \vse) } & = & \prod_{i=1}^n \exp \li [ \f{\mu - \vse}{\ba}
+ \exp \li ( \f{\vse - x_i}{\ba} \ri ) - \exp \li ( \f{\mu -
x_i}{\ba} \ri )  \ri ]\\
& = &  \li [ \exp \li ( \f{\mu - \vse}{\ba} \ri ) \ri ]^n \exp \li
\{  \li [ \exp \li ( \f{\vse}{\ba} \ri ) - \exp \li ( \f{ \mu}{\ba}
\ri )  \ri ] \sum_{i=1}^n \exp \li (  - \f{x_i}{\ba} \ri ) \ri \}.
\eee Observing that for $\vse \in (- \iy, \mu]$,
\[
\li [ \exp \li ( \f{\vse}{\ba} \ri ) - \exp \li ( \f{ \mu}{\ba} \ri
)  \ri ] \exp \li (  - \f{x}{\ba} \ri )
\]
is a concave function of $x$, we have that
\[
\f{ f_{\bs{\mcal{X}}} (\bs{x})  }{ g_{\bs{\mcal{X}}} (\bs{x}, \vse)
} \leq \li [ \exp \li ( \f{\mu - \vse}{\ba} \ri ) \ri ]^n \exp \li
\{  n \li [ \exp \li ( \f{\vse}{\ba} \ri ) - \exp \li ( \f{
\mu}{\ba} \ri )  \ri ] \exp \li (  - \f{\ovl{\bs{x}}_n }{\ba} \ri )
\ri \},
\]
where $\ovl{\bs{x}}_n = \f{\sum_{i=1}^n x_i}{n}$.  In view of the
fact that for $\vse \in (- \iy, \mu]$,
\[
\li [ \exp \li ( \f{\vse}{\ba} \ri ) - \exp \li ( \f{ \mu}{\ba} \ri
)  \ri ] \exp \li (  - \f{x}{\ba} \ri )
\]
is also an increasing function of $x$, we have that
\[
\f{ f_{\bs{\mcal{X}}} (\bs{x})  }{ g_{\bs{\mcal{X}}} (\bs{x}, \vse) } \leq \Lm (\vse) \qqu \fa \vse \in (- \iy, \mu] \; \tx{provided that} \;
\ovl{\bs{x}}_n \leq z,
\]
where
\[
\Lm (\vse) = \li [ \exp \li ( \f{\mu - \vse}{\ba} \ri ) \ri ]^n \exp
\li \{  n \li [ \exp \li ( \f{\vse}{\ba} \ri ) - \exp \li ( \f{
\mu}{\ba} \ri )  \ri ] \exp \li (  - \f{ z }{\ba} \ri ) \ri \}.
\]
This implies that $f_{\bs{\mcal{X}}} ( \bs{\mcal{X}}) \; \bb{I}_{ \{ \ovl{X}_n \leq z \} } \leq \Lm (\vse) \; g_{\bs{\mcal{X}}} ( \bs{\mcal{X}},
\vse)$ holds for any $\vse \in (- \iy, \mu]$. By virtue of Theorem \ref{ThM888}, we have $\Pr \li \{ \ovl{X}_n \leq z \ri \} \leq \inf_{ \vse
\in (- \iy, \mu] } \Lm (\vse)$.   Note that
\[
\Lm (\vse ) = \li \{ \exp \li [  w (\vse) + \f{ \mu }{\ba}  - \exp
\li ( \f{ \mu - z }{\ba} \ri )  \ri ] \ri \}^n,
\]
where
\[
w(\vse) = \f{ - \vse }{\ba} + \exp \li ( \f{\vse}{\ba} \ri )  \exp
\li (  - \f{z}{\ba} \ri ), \qu \vse \in (-\iy, \mu].
\]
It can be checked that the first and second derivatives of $w (\vse)$ are
\[
w^\prime (\vse) = \f{ - 1 }{\ba} + \f{1}{\ba} \exp \li (
\f{\vse}{\ba} \ri )  \exp \li (  - \f{z}{\ba} \ri ), \qu w^{\prime
\prime} (\vse) =  \f{1}{\ba^2} \exp \li ( \f{\vse}{\ba} \ri ) \exp
\li (  - \f{z}{\ba} \ri ).
\]
Obviously, \[ w^\prime (z) = 0, \qqu w^{\prime \prime} (z) =
\f{1}{\ba^2} > 0.
\]
It follows that
\[
\inf_{ \vse \in (- \iy, \mu] } \Lm (\vse) = \Lm (z).
\]
Therefore, \bee \Pr \{ \ovl{X}_n \leq z \} \leq \Lm (z) = \li \{ \exp \li [ \f{ \mu - z }{\ba} + 1 - \exp \li ( \f{ \mu - z }{\ba} \ri ) \ri ]
\ri \}^n \eee for $z \leq \mu$.  This completes the proof of the theorem.

\subsection{Proof of Theorem \ref{invGa}} \la{invGaapp}

Let $\bs{\mcal{X}} = [ X_1, \cd, X_n ]$ and $\bs{x} = [x_1, \cd,
x_n]$.  The joint
 probability density function of $\bs{\mcal{X}}$ is
 \[
f_{\bs{\mcal{X}}} (\bs{x}) =  \prod_{i=1}^n \f{ \ba^\al } {  \Ga (
\al) } x_i^{-\al - 1} \exp \li ( - \f{\ba}{x_i} \ri ).
 \]
To apply the LR method to show (\ref{invgamma1}), we construct a family of probability density functions
\[
g_{\bs{\mcal{X}}} (\bs{x}, \vse) = \prod_{i=1}^n \f{ \ba^\vse } {
\Ga ( \vse) } x_i^{-\vse - 1} \exp \li ( - \f{\ba}{x_i} \ri ), \qqu
\vse \in [ \al, \iy).
\]
Clearly, \bee  \f{ f_{\bs{\mcal{X}}} (\bs{x}) }{ g_{\bs{\mcal{X}}}
(\bs{x}, \vse) } & = & \prod_{i=1}^n \f{ \Ga (\vse) }{ \Ga (\al) }
\ba^{\al - \vse} x_i^{ \vse - \al  } \\
& = & \li [ \f{ \Ga (\vse) }{ \Ga (\al) } \ba^{\al - \vse} \ri ]^n
\li (  \prod_{i=1}^n x_i \ri )^{ \vse - \al }\\
& \leq & \li [ \f{ \Ga (\vse) }{ \Ga (\al) } \ba^{\al - \vse} \ri ]^n \li (  \ovl{\bs{x}}_n \ri )^{ n (\vse - \al) },  \eee where
$\ovl{\bs{x}}_n = \f{\sum_{i=1}^n x_i}{n}$.  It follows that  \[ \f{ f_{\bs{\mcal{X}}} (\bs{x})  }{ g_{\bs{\mcal{X}}} (\bs{x}, \vse) } \leq \Lm
(\vse) \qqu \fa \vse \in [ \al, \iy) \; \tx{provided that} \; \ovl{\bs{x}}_n \leq z,
\]
where
\[
\Lm (\vse) = \li [ \f{ \Ga (\vse) }{ \Ga (\al) } \li ( \f{z}{\ba}
\ri )^{\vse - \al} \ri ]^n.
\]
This implies that $f_{\bs{\mcal{X}}} ( \bs{\mcal{X}}) \; \bb{I}_{ \{
\ovl{X}_n \leq z \} } \leq \Lm (\vse) \; g_{\bs{\mcal{X}}} (
\bs{\mcal{X}}, \vse)$ holds for any $\vse \in [\al, \iy)$.  By
virtue of Theorem \ref{ThM888}, we have $\Pr \li \{ \ovl{X}_n \leq z
\ri \} \leq \inf_{ \vse \in [\al, \iy) } \Lm (\vse)$.  As a
consequence of $0 < z \leq \f{ \ba } { \al - 1 }$, we have
$\f{\ba}{z} + 1 \geq \al$.  It follows that \bee \Pr \{ \ovl{X}_n
\leq z \} \leq \Lm \li (  \f{\ba}{z} + 1 \ri )  =   \li [ \f{ \Ga
(\f{\ba}{z} + 1) }{ \Ga (\al) } \li ( \f{ z } { \beta } \ri
)^{\f{\ba}{z} - \al + 1} \ri ]^n \qu \tx{for} \; 0 < z \leq \f{ \ba
} { \al - 1 }. \eee  This proves inequality (\ref{invgamma1}).

To apply the LR method to show inequality (\ref{invgamma2}), we construct a family of probability density functions
\[
g_{\bs{\mcal{X}}} (\bs{x}, \vse) =  \prod_{i=1}^n \f{ \vse^\al } {
\Ga ( \al) } x_i^{-\al - 1} \exp \li ( - \f{\vse}{x_i} \ri ), \qqu
\vse \in (0, \ba].
\]
It can be seen that
\[
\f{ f_{\bs{\mcal{X}}} (\bs{x}) }{ g_{\bs{\mcal{X}}} (\bs{x}, \vse) }
= \li ( \f{\ba}{\vse} \ri )^{n \al}  \exp \li ( \sum_{i=1}^n
\f{\vse - \ba}{ x_i } \ri ).
\]
Observing that for $\vse \in (0, \ba]$, $\f{\vse - \ba}{ x }$ is a
concave function of $x > 0$, we have that
\[
\f{ f_{\bs{\mcal{X}}} (\bs{x}) }{ g_{\bs{\mcal{X}}} (\bs{x}, \vse) }
\leq \li [ \li ( \f{\ba}{\vse} \ri )^{ \al}  \exp \li ( \f{\vse -
\ba}{ \ovl{\bs{x}}_n } \ri ) \ri ]^n.
\]
It follows that
\[ \f{ f_{\bs{\mcal{X}}}
(\bs{x})  }{ g_{\bs{\mcal{X}}} (\bs{x}, \vse) } \leq \Lm (\vse) \qqu \fa \vse \in (0, \ba] \; \tx{provided that} \; \ovl{\bs{x}}_n \leq z,
\]
where
\[
\Lm (\vse) = \li [ \li ( \f{\ba}{\vse} \ri )^{ \al}  \exp \li (
\f{\vse - \ba}{ z } \ri ) \ri ]^n.
\]
This implies that $f_{\bs{\mcal{X}}} ( \bs{\mcal{X}}) \; \bb{I}_{ \{
\ovl{X}_n \leq z \} } \leq \Lm (\vse) \; g_{\bs{\mcal{X}}} (
\bs{\mcal{X}}, \vse)$ holds for any $\vse \in (0, \ba]$.  By virtue
of Theorem \ref{ThM888}, we have $\Pr \li \{ \ovl{X}_n \leq z \ri \}
\leq \inf_{ \vse \in (0, \ba] } \Lm (\vse)$.  By differentiation, it
can be shown that the infimum of $\Lm (\vse)$ with respect to $\vse
\in (0, \ba]$ is attained at $\vse = \al z$ as long as $0 < z \leq
\f{\ba}{\al}$.  Therefore,
\[
\Pr \li \{ \ovl{X}_n \leq z \ri \} \leq  \Lm \li (  \al z \ri ) =
\li [ \li ( \f{\ba}{\al z} \ri )^\al
 \exp \li ( \f{ \al z - \ba}{ z } \ri ) \ri ]^n  \qu \tx{for} \; 0 < z \leq \f{\ba}{\al}.
\]
This proves inequality (\ref{invgamma2}).  The proof of the theorem
is thus completed.

\subsection{Proof of Theorem \ref{ineqinvGaussian}}
\la{ineqinvGaussianapp}

Let $\bs{\mcal{X}} = [ X_1, \cd, X_n ]$ and $\bs{x} = [x_1, \cd,
x_n]$.  The joint
 probability density function of $\bs{\mcal{X}}$ is
 \[
f_{\bs{\mcal{X}}} (\bs{x}) =  \prod_{i=1}^n \li ( \f{\lm}{2 \pi
x_i^3} \ri )^{1 \sh 2}  \exp \li ( - \f{ \lm (x_i - \se)^2}{2 \se^2
x_i} \ri ).
 \]
To apply the LR method to show (\ref{inversegausineq}), we construct a family of probability density functions
\[
g_{\bs{\mcal{X}}} (\bs{x}, \vse) = \prod_{i=1}^n \li ( \f{\lm}{2 \pi
x_i^3} \ri )^{1 \sh 2}  \exp \li ( - \f{ \lm (x_i - \vse)^2}{2
\vse^2 x_i} \ri ), \qqu \vse \in ( 0, \se].
\]
It can be verified that
\[
\f{ f_{\bs{\mcal{X}}} (\bs{x}) }{ g_{\bs{\mcal{X}}} (\bs{x}, \vse) }
= \li \{ \exp \li [  \f{ \lm}{ \se}  - \f{ \lm}{ \vse} + \li ( \f{
\lm }{2 \vse^2} - \f{ \lm }{2 \se^2} \ri ) \ovl{\bs{x}}_n \ri ] \ri
\}^n,
\]
where $\ovl{\bs{x}}_n = \f{\sum_{i=1}^n x_i}{n}$.  It follows that
\[
\f{ f_{\bs{\mcal{X}}} (\bs{x}) }{ g_{\bs{\mcal{X}}} (\bs{x}, \vse) }
\leq \Lm (\vse) \qu \fa \vse \in (0, \se] \; \tx{provided that} \;
\ovl{\bs{x}}_n  \leq z,
\]
where
\[
\Lm (\vse) = \li \{ \exp \li [ \f{ \lm}{ \se}  - \f{ \lm}{ \vse} +
\li ( \f{ \lm }{2 \vse^2} - \f{ \lm }{2 \se^2} \ri ) z \ri ] \ri
\}^n.
\]
This implies that $f_{\bs{\mcal{X}}} ( \bs{\mcal{X}}) \; \bb{I}_{ \{ \ovl{X}_n \leq z \} } \leq \Lm (\vse) \; g_{\bs{\mcal{X}}} ( \bs{\mcal{X}},
\vse)$ holds for any $\vse \in (0, \se]$.  By virtue of Theorem \ref{ThM888}, we have $\Pr \li \{ \ovl{X}_n \leq z \ri \} \leq \inf_{ \vse \in
(0, \se] } \Lm (\vse)$.  By differentiation, it can be shown that the infimum of $\Lm (\vse)$ with respect to $\vse \in (0, \se]$ is attained at
$\vse = z$ as long as $0 < z \leq \se$. Therefore,
\[
\Pr \{ \ovl{X}_n \leq z \}  \leq \Lm (z) =  \li [ \exp \li (
\f{\lm}{\se} - \f{\lm}{2 z} - \f{ \lm z }{ 2 \se^2  }  \ri ) \ri ]^n
\qu \tx{for} \; 0 < z \leq \se.
\]
This completes the proof of the theorem.

\subsection{Proof of Theorem \ref{Laglog}} \la{Laglogapp}

Let $\bs{\mcal{X}} = [ X_1, \cd, X_n ]$ and $\bs{x} = [x_1, \cd,
x_n]$.  The joint
 probability mass function of $\bs{\mcal{X}}$ is
 \[
f_{\bs{\mcal{X}}} (\bs{x}) =  \prod_{i=1}^n \f{ - \se^{x_i} (1 - \se)^{x_i (\ba - 1)} \Ga (\ba x_i)  } { \Ga(x_i + 1) \Ga (\ba x_i - x_i + 1)
\ln (1 - \se) }.
 \]
To apply the LR method to show (\ref{LagBin}), we construct a family of probability mass functions
\[
g_{\bs{\mcal{X}}} (\bs{x}, \vse) = \prod_{i=1}^n \f{ - \vse^{x_i} (1 - \vse)^{x_i (\ba - 1)} \Ga (\ba x_i)  } { \Ga(x_i + 1) \Ga (\ba x_i - x_i
+ 1) \ln (1 - \vse) }, \qqu \vse \in ( 0, \se].
\]
It can be seen that
\[
\f{ f_{\bs{\mcal{X}}} (\bs{x}) }{ g_{\bs{\mcal{X}}} (\bs{x}, \vse) }
= \li [ \li (  \f{\se}{\vse}  \ri )^{\ovl{\bs{x}}_n}  \li ( \f{1 -
\se}{1 - \vse}  \ri
 )^{ (\ba - 1) \ovl{\bs{x}}_n} \f{ \ln ( 1 - \vse) } { \ln ( 1 - \se) }
\ri ]^n,
\]
where $\ovl{\bs{x}}_n = \f{\sum_{i=1}^n x_i}{n}$.  Define function
\[
h (x) = \ln x + (\ba - 1) \ln (1 - x)
\]
for $x \in (0, 1)$.   Then, we can write
\[
\f{ f_{\bs{\mcal{X}}} (\bs{x}) }{ g_{\bs{\mcal{X}}} (\bs{x}, \vse) }
= \li [ \exp \li ( [ h (\se) - h (\vse) ] \ovl{\bs{x}}_n  \ri ) \f{
\ln ( 1 - \vse) } { \ln ( 1 - \se) } \ri ]^n.
\]
Note that the first derivative of $h(x)$ is
\[
h^\prime (x) = \f{1 - \ba x}{x (1 - x)},
\]
which is positive for $x \in (0, \f{1}{\ba} )$.   Hence, \[ \f{
f_{\bs{\mcal{X}}} (\bs{x}) }{ g_{\bs{\mcal{X}}} (\bs{x}, \vse) }
\leq \Lm (\vse) \qu \fa \vse \in (0, \se] \; \tx{provided that} \;
\ovl{\bs{x}}_n  \leq z,
\]
where
\[
\Lm (\vse) = \li [ \exp \li ( [ h (\se) - h (\vse) ] z \ri ) \f{ \ln
( 1 - \vse) } { \ln ( 1 - \se) } \ri ]^n = \li [ \li (
\f{\se}{\vse}  \ri )^z  \li (  \f{1 - \se}{1 - \vse}  \ri
 )^{z (\ba - 1)} \f{ \ln ( 1 - \vse) } { \ln ( 1 - \se) }
\ri ]^n.
\]
This implies that $f_{\bs{\mcal{X}}} ( \bs{\mcal{X}}) \; \bb{I}_{ \{
\ovl{X}_n \leq z \} } \leq \Lm (\vse) \; g_{\bs{\mcal{X}}} (
\bs{\mcal{X}}, \vse)$ holds for any $\vse \in (0, \se]$.  By virtue
of Theorem \ref{ThM888}, we have $\Pr \li \{ \ovl{X}_n \leq z \ri \}
\leq \inf_{ \vse \in (0, \se] } \Lm (\vse)$.  By differentiation, it
can be shown that, as long as $0 < z \leq \f{ \se }{ (\ba \se - 1)
\ln ( 1 - \se) }$,  the infimum of $\Lm (\vse)$ with respect to
$\vse \in (0, \se]$ is attained at $\vse$ such that $z = \f{ \vse }{
(\ba \vse - 1) \ln ( 1 - \vse) }$.  Such number $\vse$ is unique
because the first derivative of $\f{ \vse }{ (\ba \vse - 1) \ln ( 1
- \vse) }$ with respective to $\vse \in (0, \f{1}{\ba} )$ is equal
to
\[ \f{ 1 }{ [ (1 - \ba \vse) \ln ( 1 - \vse) ]^2
} \li [  - \ln (1 - \vse) - \vse \f{1 - \ba \vse}{1 - \vse} \ri ],
\]
which is no less than
\[
\f{ (\ba - 1) \vse^2 }{ (1 - \vse) [ (1 - \ba
\vse) \ln ( 1 - \vse) ]^2 }  > 0.
\]
This completes the proof of the theorem.

\subsection{Proof of Theorem \ref{NegBin}} \la{NegBinapp}

Let $\bs{\mcal{X}} = [ X_1, \cd, X_n ]$ and $\bs{x} = [x_1, \cd,
x_n]$.  The joint
 probability mass function of $\bs{\mcal{X}}$ is
 \[
f_{\bs{\mcal{X}}} (\bs{x}) =  \prod_{i=1}^n \f{\ba}{ \al x_i + \ba}
\bi{\al x_i + \ba}{x_i} \se^{x_i} (1 - \se)^{\ba + \al x_i - x_i}.
 \]
To apply the LR method to show (\ref{LagNeg}), we construct a family of probability mass functions
\[
g_{\bs{\mcal{X}}} (\bs{x}, \vse) = \prod_{i=1}^n \f{\ba}{ \al x_i +
\ba} \bi{\al x_i + \ba}{x_i} \vse^{x_i} (1 - \vse)^{\ba + \al x_i -
x_i}, \qu \vse \in (0, \se ].
\]
It can be seen that
\[
\f{ f_{\bs{\mcal{X}}} (\bs{x}) }{ g_{\bs{\mcal{X}}} (\bs{x}, \vse) }
= \li [ \li (  \f{\se}{\vse}  \ri )^{\ovl{\bs{x}}_n}  \li ( \f{1 -
\se}{1 - \vse}  \ri
 )^{ \ba + (\al - 1) \ovl{\bs{x}}_n}
\ri ]^n,
\]
where $\ovl{\bs{x}}_n = \f{\sum_{i=1}^n x_i}{n}$.  Define function
\[
h (x) = \ln x + (\al - 1) \ln (1 - x)
\]
for $x \in (0, 1)$. Then, we can write
\[
\f{ f_{\bs{\mcal{X}}} (\bs{x}) }{ g_{\bs{\mcal{X}}} (\bs{x}, \vse) }
= \li [ \exp \li ( [ h (\se) - h (\vse) ] \ovl{\bs{x}}_n  \ri ) \li
( \f{1 - \se}{1 - \vse}  \ri
 )^{ \ba }  \ri ]^n.
\]
Note that the first derivative of $h(x)$ is
\[
h^\prime (x) = \f{1 - \al x}{x (1 - x)},
\]
which is positive for $x \in (0, \f{1}{\al} )$. Hence, \[ \f{
f_{\bs{\mcal{X}}} (\bs{x}) }{ g_{\bs{\mcal{X}}} (\bs{x}, \vse) }
\leq \Lm (\vse) \qu \fa \vse \in (0, \se] \; \tx{provided that} \;
\ovl{\bs{x}}_n  \leq z,
\]
where
\[
\Lm (\vse) = \li [ \exp \li ( [ h (\se) - h (\vse) ] z \ri ) \li (
\f{1 - \se}{1 - \vse}  \ri
 )^{ \ba }  \ri ]^n = \li [ \li ( \f{\se}{\vse}
\ri )^z  \li (  \f{1 - \se}{1 - \vse}  \ri
 )^{\ba + (\al - 1) z} \ri ]^n.
\]
This implies that $f_{\bs{\mcal{X}}} ( \bs{\mcal{X}}) \; \bb{I}_{ \{
\ovl{X}_n \leq z \} } \leq \Lm (\vse) \; g_{\bs{\mcal{X}}} (
\bs{\mcal{X}}, \vse)$ holds for any $\vse \in (0, \se]$.  By virtue
of Theorem \ref{ThM888}, we have $\Pr \li \{ \ovl{X}_n \leq z \ri \}
\leq \inf_{ \vse \in (0, \se] } \Lm (\vse)$.  By differentiation, it
can be shown that, as long as $0 < z \leq \f{ \ba \se  }{ 1 - \al
\se }$, the infimum of $\Lm (\vse)$ with respect to $\vse \in (0,
\se]$ is attained at $\vse = \f{z}{\ba + \al z}$. This completes the
proof of the theorem.

\subsection{Proof of Theorem \ref{laplace}} \la{laplaceapp}

Let $\bs{\mcal{X}} = [ X_1, \cd, X_n ]$ and $\bs{x} = [x_1, \cd,
x_n]$.  The joint
 probability density function of $\bs{\mcal{X}}$ is
 \[
f_{\bs{\mcal{X}}} (\bs{x}) =  \prod_{i=1}^n \f{1}{2 \ba} \exp \li (
- \f{ | x_i - \al|  }{\ba}  \ri ).
 \]
To apply the LR method to show (\ref{laplaceineqa}), we construct a family of probability density functions
\[
g_{\bs{\mcal{X}}} (\bs{x}, \vse) = \prod_{i=1}^n \f{1}{2 \vse} \exp
\li ( - \f{ | x_i - \al|  }{\vse}  \ri ), \qu \vse \in [\ba, \iy).
\]
It can be seen that for $\vse \in [\ba, \iy)$, \bee \f{
f_{\bs{\mcal{X}}} (\bs{x}) }{ g_{\bs{\mcal{X}}} (\bs{x}, \vse) } & =
& \li (  \f{\vse}{\ba}  \ri )^n \exp \li [  \li ( \f{1}{\vse} -
\f{1}{\ba} \ri )
\sum_{i=1}^n  |x_i - \al|   \ri ]\\
& \leq & \li (  \f{\vse}{\ba}  \ri )^n \exp \li [  \li ( \f{1}{\vse}
- \f{1}{\ba} \ri ) \sum_{i=1}^n  (x_i - \al) \ri ]\\
& = & \li (  \f{\vse}{\ba}  \ri )^n \li \{  \exp \li [  \li (
\f{1}{\vse} - \f{1}{\ba} \ri ) (\ovl{\bs{x}}_n - \al) \ri ] \ri
\}^n, \eee where $\ovl{\bs{x}}_n = \f{\sum_{i=1}^n x_i}{n}$.  Since
$\f{1}{\vse} - \f{1}{\ba} \leq 0$ for $\vse \in [\ba, \iy)$, it
follows that
\[ \f{
f_{\bs{\mcal{X}}} (\bs{x}) }{ g_{\bs{\mcal{X}}} (\bs{x}, \vse) }
\leq \Lm (\vse) \qu \fa \vse \in [\ba, \iy) \; \tx{provided that} \;
\ovl{\bs{x}}_n  \geq z,
\]
where
\[
\Lm (\vse) = \li \{ \f{\vse}{\ba} \exp \li [  \li ( \f{1}{\vse} -
\f{1}{\ba} \ri ) (z - \al) \ri ] \ri \}^n.
\]
This implies that $f_{\bs{\mcal{X}}} ( \bs{\mcal{X}}) \; \bb{I}_{ \{ \ovl{X}_n \geq z \} } \leq \Lm (\vse) \; g_{\bs{\mcal{X}}} ( \bs{\mcal{X}},
\vse)$ holds for any $\vse \in [\ba, \iy)$.  By virtue of Theorem \ref{ThM888}, we have $\Pr \li \{ \ovl{X}_n \geq z \ri \} \leq \inf_{ \vse \in
[\ba, \iy) } \Lm (\vse)$.  By differentiation, it can be shown that, as long as $z \geq \al + \ba $, the infimum of $\Lm (\vse)$ with respect to
$\vse \in [\ba, \iy)$ is attained at $\vse = z - \al$. This proves (\ref{laplaceineqa}).

To show (\ref{laplaceineqb}), note that for $\vse \in [\ba, \iy)$,
\bee \f{ f_{\bs{\mcal{X}}} (\bs{x}) }{ g_{\bs{\mcal{X}}} (\bs{x},
\vse) } & = & \li (  \f{\vse}{\ba}  \ri )^n \exp \li [  \li (
\f{1}{\vse} - \f{1}{\ba} \ri )
\sum_{i=1}^n  |x_i - \al|   \ri ]\\
& \leq & \li (  \f{\vse}{\ba}  \ri )^n \exp \li [  \li ( \f{1}{\vse}
- \f{1}{\ba} \ri ) \sum_{i=1}^n  (\al - x_i) \ri ]\\
& = & \li (  \f{\vse}{\ba}  \ri )^n \li \{  \exp \li [  \li (
\f{1}{\vse} - \f{1}{\ba} \ri ) (\al - \ovl{\bs{x}}_n) \ri ] \ri
\}^n. \eee Since $\f{1}{\vse} - \f{1}{\ba} \leq 0$ for $\vse \in
[\ba, \iy)$, it follows that
\[ \f{
f_{\bs{\mcal{X}}} (\bs{x}) }{ g_{\bs{\mcal{X}}} (\bs{x}, \vse) }
\leq \Lm (\vse) \qu \fa \vse \in [\ba, \iy) \; \tx{provided that} \;
\ovl{\bs{x}}_n  \leq z,
\]
where
\[
\Lm (\vse) = \li \{ \f{\vse}{\ba} \exp \li [  \li ( \f{1}{\vse} -
\f{1}{\ba} \ri ) (\al - z) \ri ] \ri \}^n.
\]
This implies that $f_{\bs{\mcal{X}}} ( \bs{\mcal{X}}) \; \bb{I}_{ \{ \ovl{X}_n \leq z \} } \leq \Lm (\vse) \; g_{\bs{\mcal{X}}} ( \bs{\mcal{X}},
\vse)$ holds for any $\vse \in [\ba, \iy)$.  By virtue of Theorem \ref{ThM888}, we have $\Pr \li \{ \ovl{X}_n \leq z \ri \} \leq \inf_{ \vse \in
[\ba, \iy) } \Lm (\vse)$.  By differentiation, it can be shown that, as long as $z \leq \al - \ba $, the infimum of $\Lm (\vse)$ with respect to
$\vse \in [\ba, \iy)$ is attained at $\vse = \al - z$. This proves (\ref{laplaceineqb}). The proof of the theorem is thus completed.

\subsection{Proof of Theorem \ref{ineqlogarithmic}}
\la{ineqlogarithmicapp}

Let $\bs{\mcal{X}} = [ X_1, \cd, X_n ]$ and $\bs{x} = [x_1, \cd,
x_n]$.  The joint
 probability mass function of $\bs{\mcal{X}}$ is
 \[
f_{\bs{\mcal{X}}} (\bs{x}) =  \prod_{i=1}^n \f{ q^{x_i}}{x_i \ln
\f{1}{1 - q}}.
 \]
To apply the LR method to show (\ref{logarithmicsingle}), we construct a family of probability mass functions
\[
g_{\bs{\mcal{X}}} (\bs{x}, \vse) = \prod_{i=1}^n \f{ \vse^{x_i}}{x_i
\ln \f{1}{1 - \vse}}, \qu \vse \in (0, q].
\]
Clearly, \bee \f{ f_{\bs{\mcal{X}}} (\bs{x}) }{ g_{\bs{\mcal{X}}}
(\bs{x}, \vse) } =  \li [ \f{ \ln (1 - q) }{ \ln (1 - \vse) } \li (
\f{q}{\vse} \ri )^{\ovl{\bs{x}}_n} \ri ]^n, \eee where
$\ovl{\bs{x}}_n = \f{\sum_{i=1}^n x_i}{n}$.  Hence,
\[ \f{
f_{\bs{\mcal{X}}} (\bs{x}) }{ g_{\bs{\mcal{X}}} (\bs{x}, \vse) }
\leq \Lm (\vse) \qu \fa \vse \in (0, q] \; \tx{provided that} \;
\ovl{\bs{x}}_n  \leq z,
\]
where
\[
\Lm (\vse) = \li [ \f{ \ln (1 - q) }{ \ln (1 - \vse) } \li (
\f{q}{\vse} \ri )^z \ri ]^n.
\]
This implies that $f_{\bs{\mcal{X}}} ( \bs{\mcal{X}}) \; \bb{I}_{ \{
\ovl{X}_n \leq z \} } \leq \Lm (\vse) \; g_{\bs{\mcal{X}}} (
\bs{\mcal{X}}, \vse)$ holds for any $\vse \in (0, q]$.  By virtue of
Theorem \ref{ThM888}, we have $\Pr \li \{ \ovl{X}_n \leq z \ri \}
\leq \inf_{ \vse \in (0, q] } \Lm (\vse)$.  By differentiation, it
can be shown that, as long as $z \leq \f{q}{(1 - q) \ln \f{1}{1 -
q}}$, the infimum of $\Lm (\vse)$ with respect to $\vse \in (0, q]$
is attained at $\vse \in (0, q]$ such that $z = \f{\vse}{(1 - \vse)
\ln \f{1}{1 - \vse}}$.   Such number $\vse$ is unique because the
function $\f{\vse}{(1 - \vse) \ln \f{1}{1 - \vse}}$ is increasing
with respect to $\vse \in (0, 1)$. The proof of the theorem is thus
completed.

\subsection{Proof of Theorem \ref{lognormal}} \la{lognormalapp}

Let $\bs{\mcal{X}} = [ X_1, \cd, X_n ]$ and $\bs{x} = [x_1, \cd,
x_n]$.  The joint
 probability density function of $\bs{\mcal{X}}$ is
 \[
f_{\bs{\mcal{X}}} (\bs{x}) =  \prod_{i=1}^n \f{1}{ x_i \sq{2 \pi}
\si }  \exp \li [  - \f{1}{2} \li ( \f{ \mu - \ln x_i }{\si} \ri )^2
\ri ].
 \]
To apply the LR method to show (\ref{lognormineq}), we construct a family of probability density functions
\[
g_{\bs{\mcal{X}}} (\bs{x}, \vse) = \prod_{i=1}^n \f{1}{ x_i \sq{2
\pi} \si }  \exp \li [  - \f{1}{2} \li ( \f{ \vse - \ln x_i }{\si}
\ri )^2 \ri ], \qu \vse \in (0, \mu].
\]
It can be seen that \bee \f{ f_{\bs{\mcal{X}}} (\bs{x}) }{
g_{\bs{\mcal{X}}} (\bs{x}, \vse) } =  \exp \li [   \f{\vse -
\mu}{\si^2} \sum_{i=1}^n \li ( \f{ \mu + \vse }{2} - \ln x_i \ri )
\ri ]. \eee It can be readily shown that for $\vse \in (0, \mu]$,
\[
\f{\vse - \mu}{\si^2}  \li ( \f{ \mu + \vse }{2} - \ln x \ri )
\]
is a concave function of $x > 0$.  Hence,
\[
\f{ f_{\bs{\mcal{X}}} (\bs{x}) }{ g_{\bs{\mcal{X}}} (\bs{x}, \vse) }
\leq  \li \{ \exp \li [   \f{\vse - \mu}{\si^2} \li ( \f{ \mu + \vse
}{2} - \ln \ovl{\bs{x}}_n \ri )^2 \ri ] \ri \}^n \qu \tx{for} \;
\vse \in (0, \mu],
\]
where $\ovl{\bs{x}}_n = \f{\sum_{i=1}^n x_i}{n}$. It follows that \[
\f{ f_{\bs{\mcal{X}}} (\bs{x}) }{ g_{\bs{\mcal{X}}} (\bs{x}, \vse) }
\leq \Lm (\vse) \qu \fa \vse \in (0, \mu] \; \tx{provided that} \;
\ovl{\bs{x}}_n  \leq z,
\]
where
\[
\Lm (\vse) = \li \{ \exp \li [   \f{\vse - \mu}{\si^2} \li ( \f{ \mu
+ \vse }{2} - \ln z \ri )^2 \ri ] \ri \}^n.
\]
This implies that $f_{\bs{\mcal{X}}} ( \bs{\mcal{X}}) \; \bb{I}_{ \{
\ovl{X}_n \leq z \} } \leq \Lm (\vse) \; g_{\bs{\mcal{X}}} (
\bs{\mcal{X}}, \vse)$ holds for any $\vse \in (0, \mu]$.  By virtue
of Theorem \ref{ThM888}, we have $\Pr \li \{ \ovl{X}_n \leq z \ri \}
\leq \inf_{ \vse \in (0, \mu] } \Lm (\vse)$.  By differentiation, it
can be shown that, as long as $0 < z \leq e^\mu$, the infimum of
$\Lm (\vse)$ with respect to $\vse \in (0, \mu]$ is attained at
$\vse = \ln z$. Therefore,
\[
\Pr \{ \ovl{X}_n \leq z \} \leq  \Lm ( \ln z ) = \exp \li [ -
\f{n}{2} \li ( \f{\mu - \ln z}{\si} \ri )^2  \ri ] \qu \tx{for $0 <
z \leq e^\mu$}.
\]
The proof of the theorem is thus completed.

\subsection{Proof of Theorem \ref{nakagami}} \la{nakagamiapp}

Let $\bs{\mcal{X}} = [ X_1, \cd, X_n ]$ and $\bs{x} = [x_1, \cd,
x_n]$.  The joint
 probability density function of $\bs{\mcal{X}}$ is
 \[
f_{\bs{\mcal{X}}} (\bs{x}) =  \prod_{i=1}^n \f{2}{\Ga(m)} \f{
x_i^{2m - 1} }{ \si^{2m} } \exp \li ( - \f{x_i^2}{\si^2} \ri ).
 \]
To apply the LR method to show (\ref{nakagamia}), we construct a family of probability density functions
\[
g_{\bs{\mcal{X}}} (\bs{x}, \vse) = \prod_{i=1}^n \f{2}{\Ga(\vse)}
\f{ x_i^{2\vse - 1} }{ \si^{2\vse} } \exp \li ( - \f{x_i^2}{\si^2}
\ri ), \qu \vse \in (0, m].
\]
Clearly, for $\vse \in (0, m]$, \bee \f{ f_{\bs{\mcal{X}}} (\bs{x})
}{ g_{\bs{\mcal{X}}} (\bs{x}, \vse) } =  \li [ \f{ \Ga ( \vse ) }{
\Ga (m) } \si^{2 (\vse - m)} \ri ]^n \li ( \prod_{i=1}^n x_i \ri )^{
2 (m - \vse)  } \leq \li [ \f{ \Ga ( \vse ) }{ \Ga (m) } \si^{2
(\vse - m)} \ri ]^n \li ( \ovl{\bs{x}}_n \ri )^{ 2 n (m - \vse)  },
\eee where $\ovl{\bs{x}}_n = \f{\sum_{i=1}^n x_i}{n}$.  It follows
that \[ \f{ f_{\bs{\mcal{X}}} (\bs{x}) }{ g_{\bs{\mcal{X}}} (\bs{x},
\vse) } \leq \Lm (\vse) \qu \fa \vse \in (0, m] \; \tx{provided
that} \; \ovl{\bs{x}}_n  \leq z,
\]
where
\[
\Lm (\vse) =  \li [ \f{ \Ga ( \vse ) }{ \Ga (m) } \li ( \f{z}{\si}
\ri )^{2 (m - \vse)} \ri ]^n.
\]
This implies that $f_{\bs{\mcal{X}}} ( \bs{\mcal{X}}) \; \bb{I}_{ \{
\ovl{X}_n \leq z \} } \leq \Lm (\vse) \; g_{\bs{\mcal{X}}} (
\bs{\mcal{X}}, \vse)$ holds for any $\vse \in (0, m]$.  By virtue of
Theorem \ref{ThM888}, we have $\Pr \li \{ \ovl{X}_n \leq z \ri \}
\leq \inf_{ \vse \in (0, \mu] } \Lm (\vse)$.  Letting $z = \f{ \Ga (
\vse + \f{1}{2} ) }{ \Ga (\vse) } \si $ leads to (\ref{nakagamia}).

To apply the LR method to show (\ref{nakagamib}), we construct a family of probability density functions
\[
g_{\bs{\mcal{X}}} (\bs{x}, \vse) =  \prod_{i=1}^n \f{2}{\Ga(m)} \f{
x_i^{2m - 1} }{ \vse^{2m} } \exp \li ( - \f{x_i^2}{\vse^2} \ri ),
\qu \vse \in [\si, \iy).
\]
It can be seen that
\[
\f{ f_{\bs{\mcal{X}}} (\bs{x}) }{ g_{\bs{\mcal{X}}} (\bs{x}, \vse) }
= \li ( \f{\vse}{\si} \ri )^{2 m n} \exp \li [ \li ( \f{1}{\vse^2} -
\f{1}{\si^2} \ri ) \sum_{i=1}^n  x_i^2 \ri ].
\]
Observing that for $\vse \in [\si, \iy)$, $\li ( \f{1}{\vse^2} -
\f{1}{\si^2} \ri ) x^2$ is a concave function of $x > 0$, we have
that \[ \f{ f_{\bs{\mcal{X}}} (\bs{x}) }{ g_{\bs{\mcal{X}}} (\bs{x},
\vse) } \leq \li \{  \li ( \f{\vse}{\si} \ri )^{2 m } \exp \li [ \li
( \f{1}{\vse^2} - \f{1}{\si^2} \ri ) ( \ovl{\bs{x}}_n )^2 \ri ] \ri
\}^n, \qu \fa \vse \in [\si, \iy).
\]
It follows that \[ \f{ f_{\bs{\mcal{X}}} (\bs{x}) }{
g_{\bs{\mcal{X}}} (\bs{x}, \vse) } \leq \Lm (\vse) \qu \fa \vse \in
[\si, \iy) \; \tx{provided that} \; \ovl{\bs{x}}_n  \geq z,
\]
where
\[
\Lm (\vse) =  \li [ \li ( \f{\vse}{\si} \ri )^{2 m} \exp \li (
\f{z^2}{\vse^2}  - \f{z^2}{\si^2} \ri ) \ri ]^n.
\]
This implies that $f_{\bs{\mcal{X}}} ( \bs{\mcal{X}}) \; \bb{I}_{ \{ \ovl{X}_n \geq z \} } \leq \Lm (\vse) \; g_{\bs{\mcal{X}}} ( \bs{\mcal{X}},
\vse)$ holds for any $\vse \in [\si, \iy)$.  By virtue of Theorem \ref{ThM888}, we have $\Pr \li \{ \ovl{X}_n \geq z \ri \} \leq \inf_{ \vse \in
[\si, \iy) } \Lm (\vse)$.  By differentiation, it can be shown that, as long as $z \geq \sq{m} \si$, the infimum of $\Lm (\vse)$ with respect to
$\vse \in [\si, \iy)$ is attained at $\vse = \f{z}{\sq{m}}$. Therefore,
\[
\Pr \{ \ovl{X}_n \geq z \} \leq \Lm \li ( \f{z}{\sq{m}} \ri ) = \li [ \li ( \f{z^2}{m \si^2} \ri )^{m} \exp \li ( m  - \f{z^2}{\si^2} \ri ) \ri
]^n \qu \tx{for $z \geq \sq{m} \si$}.
\]
This establishes (\ref{nakagamib}) and completes the proof of the
theorem.

\subsection{Proof of Theorem \ref{pareto}} \la{paretoapp}

Let $\bs{\mcal{X}} = [ X_1, \cd, X_n ]$ and $\bs{x} = [x_1, \cd,
x_n]$.  The joint
 probability density function of $\bs{\mcal{X}}$ is
 \[
f_{\bs{\mcal{X}}} (\bs{x}) =  \prod_{i=1}^n \f{\se}{a}  \li (
\f{a}{x_i} \ri )^{\se + 1}.
 \]
To apply the LR method to show (\ref{paretoineq}), we construct a family of probability density functions
\[
g_{\bs{\mcal{X}}} (\bs{x}, \vse) = \prod_{i=1}^n \f{\vse}{a}  \li (
\f{a}{x_i} \ri )^{\vse + 1}, \qu \vse \in [\se,  \iy).
\]
Clearly,
\[
\f{ f_{\bs{\mcal{X}}} (\bs{x}) }{ g_{\bs{\mcal{X}}} (\bs{x}, \vse) }
= \li (  \f{ \se } {  \vse  } \ri )^n \li ( \prod_{i=1}^n x_i \ri
)^{\vse - \se} \leq  \li [  \f{ \se } {  \vse  } \li ( \f{
\ovl{\bs{x}}_n }{a} \ri )^{\vse - \se} \ri ]^n,
\]
where $\ovl{\bs{x}}_n = \f{\sum_{i=1}^n x_i}{n}$.  It follows that
\[ \f{ f_{\bs{\mcal{X}}} (\bs{x}) }{ g_{\bs{\mcal{X}}} (\bs{x},
\vse) } \leq \Lm (\vse) \qu \fa \vse \in [\se, \iy) \; \tx{provided
that} \; \ovl{\bs{x}}_n  \leq z,
\]
where
\[
\Lm (\vse) =   \li [  \f{ \se } {  \vse  } \li ( \f{z}{a} \ri
)^{\vse - \se} \ri ]^n.
\]
This implies that $f_{\bs{\mcal{X}}} ( \bs{\mcal{X}}) \; \bb{I}_{ \{
\ovl{X}_n \leq z \} } \leq \Lm (\vse) \; g_{\bs{\mcal{X}}} (
\bs{\mcal{X}}, \vse)$ holds for any $\vse \in [\se, \iy)$.  By
virtue of Theorem \ref{ThM888}, we have $\Pr \li \{ \ovl{X}_n \leq z
\ri \} \leq \inf_{ \vse \in [\se, \iy) } \Lm (\vse)$.  Hence,  for
$\ga > 1$, \[ \Pr \{ \ovl{X}_n \leq \ga a \} \leq \inf_{\vse \geq
\se} \li ( \f{ \se } { \vse  }  \ga^{\vse - \se}  \ri )^n = \se^n \;
\inf_{\vse \geq \se} \exp [ n \; w (\vse) ],
\]
where $w (\vse) = - \ln  \vse   + ( \vse - \se ) \ln  \ga$.

Now consider the minimization of $w (\vse)$ subject to $\vse \geq
\se$. Note that the first and second derivatives of $w(\vse)$ are
$w^\prime (\vse) = - \f{1}{\vse} + \ln \ga$ and $w^{\prime \prime}
(\vse) = \f{1}{\vse^2}$, respectively. Hence, the minimum is
achieved at $\vse^* = \f{1}{\ln \ga}$ provided that $1 < \ga \leq
e^{1 \sh \se}$. Accordingly, $w (\vse^*) = 1 + \ln \ln \ga - \se \ln
\ga$ and
\[
\Pr \{ \ovl{X}_n \leq \ga a \} \leq \li ( \f{ e \se }{\ga^\se} \ln
\ga \ri )^n \qu \tx{for $1 < \ga \leq e^{1 \sh \se}$}.
\]
Note that the mean of $X$ is $\mu = \f{\se a}{\se - 1}$. Letting $\ga  = \f{\ro \mu}{a}$ yields
\[
\Pr \{ \ovl{X}_n \leq \ro \mu \} \leq \li [ e \se \li ( \f{\se - 1}{\ro \se} \ri )^\se \ln \li ( \f{\ro \se}{\se - 1}  \ri ) \ri ]^n \qu \tx{for
$1 - \f{1}{\se} < \ro \leq \li ( 1 - \f{1}{\se} \ri ) \exp ( \f{1}{\se} )$}.
\]
This establishes (\ref{paretoineq}) and completes the proof of the theorem.

 \subsection{Proof of Theorem \ref{powerlaw}} \la{powerlawapp}

 Let $\bs{\mcal{X}} = [ X_1, \cd, X_n ]$ and $\bs{x} = [x_1, \cd,
x_n]$.  The joint
 probability density function of $\bs{\mcal{X}}$ is
 \[
f_{\bs{\mcal{X}}} (\bs{x}) =  \prod_{i=1}^n \f{x_i^{-\al}}{C(\al)}.
 \]
To apply the LR method to show (\ref{powineq}), we construct a probability density functions
\[
g_{\bs{\mcal{X}}} (\bs{x}, \vse) = \prod_{i=1}^n
\f{x_i^{-\vse}}{C(\vse)}, \qu \vse \in [\al,  \iy).
\]
Clearly,
\[
\f{ f_{\bs{\mcal{X}}} (\bs{x}) }{ g_{\bs{\mcal{X}}} (\bs{x}, \vse) }
= \li [ \f{ C ( \vse ) } { C(\al)  } \ri ]^n \li ( \prod_{i=1}^n x_i
\ri )^{\vse - \al} \leq  \li [ \f{ C ( \vse ) } { C(\al)  } \ri ]^n
 \li [  \li ( \ovl{\bs{x}}_n  \ri )^{\vse - \al} \ri ]^n \leq \Lm (\vse)
\]
provided that $\ovl{\bs{x}}_n \leq z$, where $\ovl{\bs{x}}_n =
\f{\sum_{i=1}^n x_i}{n}$ and \[ \Lm (\vse) = \li [ \f{ C ( \vse ) }
{ C(\al)  } z^{\vse - \al} \ri ]^n.
\]  This implies that $f_{\bs{\mcal{X}}}
( \bs{\mcal{X}}) \; \bb{I}_{ \{ \ovl{X}_n \leq z \} } \leq \Lm
(\vse) \; g_{\bs{\mcal{X}}} ( \bs{\mcal{X}}, \vse)$ holds for any
$\vse \in [\al, \iy)$.  By virtue of Theorem \ref{ThM888}, we have
$\Pr \li \{ \ovl{X}_n \leq z \ri \} \leq \Lm (\vse)$. This completes
the proof of the theorem.

\subsection{Proof of Theorem \ref{Stirling}} \la{Stirlingapp}

Let $\bs{\mcal{X}} = [ X_1, \cd, X_n ]$ and $\bs{x} = [x_1, \cd,
x_n]$.  The joint
 probability mass function of $\bs{\mcal{X}}$ is
 \[
f_{\bs{\mcal{X}}} (\bs{x}) =  \prod_{i=1}^n \f{  m! | s(x_i, m)|
\se^{x_i} }{  x_i ! [- \ln (1 - \se) ]^m }.
 \]
To apply the LR method to show (\ref{Stringineq}), we construct a family of probability mass functions
\[
g_{\bs{\mcal{X}}} (\bs{x}, \vse) = \prod_{i=1}^n \f{  m! | s(x_i, m)| \vse^{x_i} }{  x_i ! [- \ln (1 - \vse) ]^m }, \qu \vse \in (0, \se ].
\]
Clearly,
\[
\f{ f_{\bs{\mcal{X}}} (\bs{x}) }{ g_{\bs{\mcal{X}}} (\bs{x}, \vse) }
= \li [  \f{ \ln (1 - \vse) } {  \ln ( 1 - \se)  } \ri ]^{n m}  \li
[ \li ( \f{\se}{\vse} \ri )^{ \ovl{\bs{x}}_n }  \ri ]^n,
\]
where $\ovl{\bs{x}}_n = \f{\sum_{i=1}^n x_i}{n}$.  It follows that
\[ \f{ f_{\bs{\mcal{X}}} (\bs{x}) }{ g_{\bs{\mcal{X}}} (\bs{x},
\vse) } \leq \Lm (\vse) \qu \fa \vse \in (0, \se] \; \tx{provided
that} \; \ovl{\bs{x}}_n  \leq z,
\]
where
\[
\Lm (\vse) =   \li [  \f{ \ln (1 - \vse) } {  \ln ( 1 - \se)  } \ri
]^{n m}  \li [ \li ( \f{\se}{\vse} \ri )^z  \ri ]^n.
\]
This implies that $f_{\bs{\mcal{X}}} ( \bs{\mcal{X}}) \; \bb{I}_{ \{
\ovl{X}_n \leq z \} } \leq \Lm (\vse) \; g_{\bs{\mcal{X}}} (
\bs{\mcal{X}}, \vse)$ holds for any $\vse \in (0, \se]$.  By virtue
of Theorem \ref{ThM888}, we have $\Pr \li \{ \ovl{X}_n \leq z \ri \}
\leq \inf_{ \vse \in (0, \se] } \Lm (\vse)$.  By differentiation, it
can be shown that, as long as $z \leq \f{m \se} {(\se - 1) \ln (1 -
\se)}$, the infimum of $\Lm (\vse)$ with respective to $\vse \in (0,
\se]$ is attained at a number $\vse$ such that $z = \f{m \vse}
{(\vse - 1) \ln (1 - \vse)}$.  Such number $\vse$ is unique because
$\f{\vse} {(\vse - 1) \ln (1 - \vse)}$ is an increasing function of
$\vse \in (0, 1)$.   This completes the proof of the theorem.

\subsection{Proof of Theorem \ref{Snedecor}} \la{Snedecorapp}

To apply the LR method, we introduce a family of probability density
functions \[ g (x, \vse) = \f{1}{\vse}  f \li ( \f{x}{\vse} \ri ),
\qqu \vse > 0.
\]
Clearly,
\[
\f{f (x)}{g (x, \vse)} = \f{   f(x) } { \f{1}{\vse} f (\f{x}{\vse} )
} = \vse^{m \sh 2}  \li (  \f{ n + \f{m x}{\vse} }{ n + m x }  \ri
)^{(n + m) \sh 2} = \vse^{m \sh 2}  \li ( 1 +  \f{ \f{1}{\vse} - 1}{
1 + \f{n}{m x} }  \ri )^{(n + m) \sh 2}.
\]
To show inequality (\ref{Snedecorineq88a}), note that $\f{f (x)}{g
(x, \vse)}$ is decreasing with respect to $x > 0$ for $\vse \geq 1$.
Hence,
\[
\f{f (x)}{g (x, \vse)} \leq \Lm(\vse) \qu \fa \vse \in [1, \iy) \;
\tx{provided that} \; x \geq z,
\]
where \be \la{Lm888}
 \Lm (\vse) = \vse^{m \sh 2}  \li ( 1 +  \f{
\f{1}{\vse} - 1}{ 1 + \f{n}{m z} }  \ri )^{(n + m) \sh 2}. \ee
 This implies that $f (X) \; \bb{I}_{ \{ X \geq z \} } \leq \Lm (\vse) \;
g ( X, \vse)$ holds for any $\vse \in [1, \iy)$. By virtue of
Theorem \ref{ThM888}, we have
\[
\Pr \li \{ X \geq z \ri \} \leq \inf_{ \vse \in [1, \iy) } \Lm (\vse) = \Lm (z) = z^{m \sh 2}  \li (  \f{ n + m }{ n + m z } \ri )^{(n + m) \sh
2} \qqu \tx{for} \; z \geq 1.
\]

To show inequality (\ref{Snedecorineq88b}), note that  $\f{f (x)}{g
(x, \vse)}$ is increasing with respect to $x > 0$ for $0 < \vse \leq
1$. Hence,
\[
\f{f (x)}{g (x, \vse)} \leq \Lm(\vse) \qu \fa \vse \in (0, 1 ] \; \tx{provided that} \; x \leq z,
\]
where $\Lm (\vse)$ is defined by (\ref{Lm888}). This implies that $f (X) \; \bb{I}_{ \{ X \leq z \} } \leq \Lm (\vse) \; g ( X, \vse)$ holds for
any $\vse \in (0, 1]$. By virtue of Theorem \ref{ThM888}, we have
\[
\Pr \li \{ X \leq z \ri \} \leq \inf_{ \vse \in (0, 1] } \Lm (\vse) = \Lm (z) = z^{m \sh 2}  \li (  \f{ n + m }{ n + m z } \ri )^{(n + m) \sh 2}
\qqu \tx{for} \; 0 < z \leq 1.
\]
This proves inequality (\ref{Snedecorineq88b}) and completes the
proof of the theorem.

\subsection{Proof of Theorem \ref{CDFineq}} \la{CDFineqapp}

To apply the LR method, we introduce a family of probability density
functions \[ g (x, \vse) = \f{1}{\vse}  f \li ( \f{x}{\vse} \ri ),
\qqu \vse > 0.
\]
Clearly,
\[
\f{f (x)}{g (x, \vse)} = \f{   f(x) } { \f{1}{\vse} f (\f{x}{\vse} )
} = \vse \li [  \f{ n + ( \f{x}{\vse}  )^2}{ n + x^2 }  \ri ]^{(n +
1) \sh 2} = \vse \li ( 1 + \f{ \f{1}{ \vse^2} - 1 } { \f{n}{x^2} + 1
} \ri )^{ (n + 1) \sh 2 }.
\]
To show inequality (\ref{stuineq88a}), note that $\f{f (x)}{g (x,
\vse)}$ is decreasing with respect to $|x|$ for $\vse \geq 1$.
Hence,
\[
\f{f (x)}{g (x, \vse)} \leq \Lm(\vse) \qu \fa \vse \in [1, \iy) \;
\tx{provided that} \; |x| \geq z,
\]
where \be \la{Lm8889}
 \Lm (\vse) = \vse \li ( 1 + \f{ \f{1}{ \vse^2} - 1 } { \f{n}{z^2} + 1
} \ri )^{ (n + 1) \sh 2 }. \ee
 This implies that $f (X) \; \bb{I}_{ \{ |X| \geq z \} } \leq \Lm (\vse) \;
g ( X, \vse)$ holds for any $\vse \in [1, \iy)$. By virtue of
Theorem \ref{ThM888}, we have
\[
\Pr \li \{ X \geq z \ri \} \leq \inf_{ \vse \in [1, \iy) } \Lm (\vse) = \Lm (z) = z \li (  \f{ n + 1}{ n + z^2 }  \ri )^{(n + 1) \sh 2} \qqu
\tx{for} \; z \geq 1.
\]
This proves inequality (\ref{stuineq88a}).

To show inequality (\ref{stuineq88b}), note that $\f{f (x)}{g (x,
\vse)}$ is increasing with respect to $|x|$ for $\vse \in (0, 1]$.
Hence,
\[
\f{f (x)}{g (x, \vse)} \leq \Lm(\vse) \qu \fa \vse \in (0, 1]  \;
\tx{provided that} \; |x| \leq z,
\]
where $\Lm (\vse)$ is defined by (\ref{Lm8889}).  This implies that $f (X) \; \bb{I}_{ \{ |X| \leq z \} } \leq \Lm (\vse) \; g ( X, \vse)$ holds
for any $\vse \in (0, 1 ]$. By virtue of Theorem \ref{ThM888}, we have
\[
\Pr \li \{ X \leq z \ri \} \leq \inf_{ \vse \in (0, 1 ] } \Lm (\vse) = \Lm (z) = z \li (  \f{ n + 1}{ n + z^2 }  \ri )^{(n + 1) \sh 2} \qqu
\tx{for} \; 0 < z \leq 1.
\]
This proves inequality (\ref{stuineq88b}) and completes the proof of
the theorem.

\subsection{Proof of Theorem \ref{truncateexp}} \la{truncateexpapp}

Let $\bs{\mcal{X}} = [ X_1, \cd, X_n ]$ and $\bs{x} = [x_1, \cd,
x_n]$.  The joint
 probability density function of $\bs{\mcal{X}}$ is
 \[
f_{\bs{\mcal{X}}} (\bs{x}) =  \prod_{i=1}^n \f{\se}{e^\se - 1}
e^{\se x_i}.
 \]
To apply the LR method to show (\ref{expch}), we construct a family of probability density functions
\[
g_{\bs{\mcal{X}}} (\bs{x}, \vse) = \prod_{i=1}^n \f{\vse}{e^\vse -
1} e^{\vse x_i}, \qu \vse \in (-\iy, \se], \qu \vse \neq 0.
\]
Clearly,
\[
\f{ f_{\bs{\mcal{X}}} (\bs{x}) }{ g_{\bs{\mcal{X}}} (\bs{x}, \vse) }
= \li (  \f{\se}{\vse} \f{ e^\vse - 1 }{ e^\se - 1 }  \ri )^n  \exp
\li [ n (\se - \vse) \ovl{\bs{x}}_n \ri ],
\]
where $\ovl{\bs{x}}_n = \f{\sum_{i=1}^n x_i}{n}$.  It follows that
\[ \f{ f_{\bs{\mcal{X}}} (\bs{x}) }{ g_{\bs{\mcal{X}}} (\bs{x},
\vse) } \leq \Lm (\vse) \qu \fa \vse \in (- \iy, \se], \; \vse \neq 0 \; \tx{provided that} \; \ovl{\bs{x}}_n  \leq z,
\]
where
\[
\Lm (\vse) = \li (  \f{\se}{\vse} \f{ e^\vse - 1 }{ e^\se - 1 }  \ri
)^n   \exp \li [ n (\se - \vse) z \ri ].
\]
This implies that $f_{\bs{\mcal{X}}} ( \bs{\mcal{X}}) \; \bb{I}_{ \{ \ovl{X}_n \leq z \} } \leq \Lm (\vse) \; g_{\bs{\mcal{X}}} ( \bs{\mcal{X}},
\vse)$ holds for any $\vse \in (- \iy, \se], \; \vse \neq 0$. By virtue of Theorem \ref{ThM888}, we have $\Pr \li \{ \ovl{X}_n \leq z \ri \}
\leq \inf_{ \vse \in (- \iy, \se] } \Lm (\vse)$. By differentiation, it can be shown that, as long as $0 < z \leq 1 + \f{1}{e^\se - 1} -
\f{1}{\se}$ and $z \neq \f{1}{2}$, the infimum of $\Lm (\vse)$ with respective to $\vse \in (- \iy, \se], \; \vse \neq 0$ is attained at a
number $\vse \in (-\iy, \se], \; \vse \neq 0$ such that $z = 1 + \f{1}{e^\vse - 1} - \f{1}{\vse}$. Such a number is unique because
\[ \lim_{\vse \to - \iy} \li ( 1 + \f{1}{e^\vse - 1} - \f{1}{\vse}
\ri ) = 0
\]
and $1 + \f{1}{e^\vse - 1} - \f{1}{\vse}$ is increasing with respect
to $\vse \neq 0$.  To show such monotonicity, note that the first
derivative of $1 + \f{1}{e^\vse - 1} - \f{1}{\vse}$ with respective
to $\vse$ is equal to
\[ \li [ \f{ e^{\vse \sh 2} - e^{-\vse \sh 2} - \vse}{\vse (e^{\vse
\sh 2} - e^{- \vse \sh 2})} \ri ]  \li ( \f{1}{\vse} + \f{ e^{\vse
\sh 2} }{ e^\vse - 1 } \ri ),
\]
where $e^{\vse \sh 2} - e^{-\vse \sh 2} - \vse$ is a function of $\vse$ with its first derivative assuming value $0$ at $\vse = 0$, and its
second derivative equal to $ \f{1}{4} ( e^{\vse \sh 2} - e^{- \vse \sh 2} )$.   This establishes inequality (\ref{expch}).  To show
(\ref{expch88b}), it suffices to note that as $z \to \f{1}{2}$, the root of equation $z = 1 + \f{1}{e^\vse - 1} - \f{1}{\vse}$ with respect to
$\vse$ tends to $0$. This completes the proof of the theorem.

\subsection{Proof of Theorem \ref{UniformCDF}} \la{UniformCDFapp}

Let $\bs{\mcal{X}} = [ X_1, \cd, X_n ]$ and $\bs{x} = [x_1, \cd,
x_n]$.  The joint
 probability density function of $\bs{\mcal{X}}$ is $f_{\bs{\mcal{X}}} (\bs{x}) =
 1$.
To apply the LR method to show (\ref{uniformineqa}), we construct a family of probability density functions
\[
g_{\bs{\mcal{X}}} (\bs{x}, \vse) = \prod_{i=1}^n \f{\vse}{e^\vse -
1} e^{\vse x_i}, \qu \vse > 0.
\]
Clearly,
\[
\f{ f_{\bs{\mcal{X}}} (\bs{x}) }{ g_{\bs{\mcal{X}}} (\bs{x}, \vse) }
= \li [  \f{ e^\vse - 1 }{ \vse}  \exp ( - \vse \; \ovl{\bs{x}}_n )
\ri ]^n,
\]
where $\ovl{\bs{x}}_n = \f{\sum_{i=1}^n x_i}{n}$.  It follows that
\[ \f{ f_{\bs{\mcal{X}}} (\bs{x}) }{ g_{\bs{\mcal{X}}} (\bs{x},
\vse) } \leq \Lm (\vse) \qu \fa \vse > 0 \; \tx{provided that} \;
\ovl{\bs{x}}_n  \geq z,
\]
where
\[
\Lm (\vse) = \li (  \f{ e^\vse - 1 }{ \vse e^{\vse z} }  \ri )^n.
\]
This implies that $f_{\bs{\mcal{X}}} ( \bs{\mcal{X}}) \; \bb{I}_{ \{
\ovl{X}_n \geq z \} } \leq \Lm (\vse) \; g_{\bs{\mcal{X}}} (
\bs{\mcal{X}}, \vse)$ holds for any $\vse > 0$.  By virtue of
Theorem \ref{ThM888}, we have $\Pr \li \{ \ovl{X}_n \geq z \ri \}
\leq \inf_{ \vse > 0 } \Lm (\vse)$.  By differentiation, it can be
shown that, as long as $1 > z \geq \f{1}{2}$, the infimum of $\Lm
(\vse)$ with respective to $\vse > 0$ is attained at a positive
number $\vse^*$ such that $z = 1 + \f{1}{e^{\vse^*} - 1} -
\f{1}{\vse^*}$. Such a number is unique because
\[ \lim_{\vse \downarrow 0} \li ( 1 + \f{1}{e^\vse - 1} - \f{1}{\vse}
\ri ) = \f{1}{2}
\]
and $1 + \f{1}{e^\vse - 1} - \f{1}{\vse}$ is increasing with respect
to $\vse > 0$.  Therefore, we have shown that
\[ \Pr \{ \ovl{X}_n \geq z \} \leq \Lm (\vse^*) \qu \tx{for} \;
\f{1}{2} < z < 1.
\]
On the other hand, it can be shown that
\[
\Lm (\vse^*) = \li (  \inf_{s > 0} e^{-z s} \bb{E} [ e^{s X} ] \ri
)^n
\]
To establish an upper bound on $\Lm (\vse^*)$, we can use the following inequality  due to Chen \cite[Appendix H]{ChenMag},
\[
\bb{E} [ e^{s X} ] < \exp \li (  \f{s^2}{24} + \f{s}{2} \ri ), \qqu \fa s \in (-\iy, \iy).
\]
By differentiation, it can be shown that
\[
\Pr \{ \ovl{X} \geq z \} \leq \Lm (\vse^*) \leq \li [ \inf_{s > 0}
\exp \li ( \f{s^2}{24} + \f{s}{2} - z s \ri ) \ri ]^n \leq \exp \li
( - 6 n \li ( z - \f{1}{2} \ri )^2 \ri ), \qqu 1>  z > \f{1}{2}.
\]
This establishes (\ref{uniformineqa}). By a similar argument, we can show (\ref{uniformineqb}).  This completes the proof of the theorem.

\subsection{Proof of Theorem \ref{WeibullCDF}} \la{WeibullCDFapp}

Let $\bs{\mcal{X}} = [ X_1, \cd, X_n ]$ and $\bs{x} = [x_1, \cd,
x_n]$.  The joint
 probability density function of $\bs{\mcal{X}}$ is
 \[
f_{\bs{\mcal{X}}} (\bs{x}) =  \prod_{i=1}^n \al \ba x_i^{\ba - 1}
\exp \li ( - \al x_i^\ba \ri ).
 \]
To apply the LR method, we construct a family of probability density functions
\[
g_{\bs{\mcal{X}}} (\bs{x}, \vse) = \prod_{i=1}^n \vse \ba x_i^{\ba -
1} \exp \li ( - \vse x_i^\ba \ri ), \qu \vse \in (0,  \iy).
\]
Clearly,
\[
\f{ f_{\bs{\mcal{X}}} (\bs{x}) }{ g_{\bs{\mcal{X}}} (\bs{x}, \vse) }
= \li ( \f{\al}{\vse} \ri )^n \exp \li [ (\vse - \al) \sum_{i=1}^n
x_i^\ba \ri ].
\]
To show inequality (\ref{Weib88a}) under the condition that $\al
z^\ba \leq 1$ and $0 < \ba \leq 1$, we restrict $\vse$ to be no less
than $\al$.  As a consequence of $0 < \ba \leq 1$ and $\vse \geq
\al$, we have that $(\vse - \al) x^\ba$ is a concave function of $x
> 0$.  By virtue of such concavity, we have
\[
\f{ f_{\bs{\mcal{X}}} (\bs{x}) }{ g_{\bs{\mcal{X}}} (\bs{x}, \vse) }
\leq \li \{ \f{\al}{\vse} \exp \li [ (\vse - \al) ( \ovl{\bs{x}}_n
)^\ba \ri ] \ri \}^n , \qu \fa \vse \in [\al, \iy),
\]
where $\ovl{\bs{x}}_n = \f{\sum_{i=1}^n x_i}{n}$. It follows that
\[ \f{ f_{\bs{\mcal{X}}} (\bs{x}) }{ g_{\bs{\mcal{X}}} (\bs{x},
\vse) } \leq \Lm (\vse) \qu \fa \vse \in [\al, \iy) \; \tx{provided
that} \; \ovl{\bs{x}}_n  \leq z,
\]
where \be \la{reuseLm}
 \Lm (\vse) =   \li \{ \f{\al}{\vse} \exp \li
[ (\vse - \al) z^\ba \ri ] \ri \}^n. \ee This implies that
$f_{\bs{\mcal{X}}} ( \bs{\mcal{X}}) \; \bb{I}_{ \{ \ovl{X}_n \leq z
\} } \leq \Lm (\vse) \; g_{\bs{\mcal{X}}} ( \bs{\mcal{X}}, \vse)$
holds for any $\vse \in [\al, \iy)$.  By virtue of Theorem
\ref{ThM888}, we have $\Pr \li \{ \ovl{X}_n \leq z \ri \} \leq
\inf_{ \vse \in [\al, \iy) } \Lm (\vse)$.  By differentiation, it
can be shown that, as long as $\al z^\ba \leq 1$, the infimum of
$\Lm (\vse)$ with respect to $\vse \in [\al, \iy)$ is attained at
$\vse = z^{- \ba}$.  Therefore,
\[
\Pr \{ \ovl{X}_n \leq z \} \leq  \Lm ( z^{- \ba} ) = \li [  \al z^{\ba} \exp (1 - \al z^\ba ) \ri ]^n \qu \tx{for $\al z^\ba \leq 1$ and $\ba <
1$}.
\]
This proves inequality (\ref{Weib88a}).

To show inequality (\ref{Weib88b}) under the condition that $\al z^\ba \geq 1$ and $\ba > 1$, we restrict $\vse$ to be a positive number less
than $\al$.  As a consequence of $\ba > 1$ and $0 < \vse < \al$, we have that $(\vse - \al) x^\ba$ is a concave function of $x
> 0$.  By virtue of such concavity, we have
\[
\f{ f_{\bs{\mcal{X}}} (\bs{x}) }{ g_{\bs{\mcal{X}}} (\bs{x}, \vse) }
\leq \li \{ \f{\al}{\vse} \exp \li [ (\vse - \al) ( \ovl{\bs{x}}_n
)^\ba \ri ] \ri \}^n , \qu \fa \vse \in (0, \al).
\]
It follows that
\[ \f{ f_{\bs{\mcal{X}}} (\bs{x}) }{ g_{\bs{\mcal{X}}} (\bs{x},
\vse) } \leq \Lm (\vse) \qu \fa \vse \in (0, \al) \; \tx{provided
that} \; \ovl{\bs{x}}_n  \geq z,
\]
where $\Lm (\vse)$ is defined by (\ref{reuseLm}).  This implies that
$f_{\bs{\mcal{X}}} ( \bs{\mcal{X}}) \; \bb{I}_{ \{ \ovl{X}_n \geq z
\} } \leq \Lm (\vse) \; g_{\bs{\mcal{X}}} ( \bs{\mcal{X}}, \vse)$
holds for any $\vse \in (0, \al)$.  By virtue of Theorem
\ref{ThM888}, we have $\Pr \li \{ \ovl{X}_n \geq z \ri \} \leq
\inf_{ \vse \in (0, \al) } \Lm (\vse)$.  By differentiation, it can
be shown that, as long as $\al z^\ba \geq 1$, the infimum of $\Lm
(\vse)$ with respect to $\vse \in (0, \al)$ is attained at $\vse =
z^{- \ba}$.  Therefore,
\[
\Pr \{ \ovl{X}_n \geq z \} \leq  \Lm ( z^{- \ba} ) = \li [  \al
z^{\ba} \exp (1 - \al z^\ba ) \ri ]^n \qu \tx{for $\al z^\ba \geq 1$
and $\ba > 1$}.
\]
This proves inequality (\ref{Weib88b}).  The proof of the theorem is thus completed.

\section{Proofs of Multivariate Inequalities}

\subsection{Proof of Theorem  \ref{DeriCompound}}
\la{DeriCompoundapp}

To apply the LR method to show (\ref{DeC88}),  we introduce a family of probability mass functions
\[
g(\bs{x}, \bs{\vse}) = \bi{n}{ \bs{x} } \f{ \Ga (\sum_{\ell = 0}^k
\vse_\ell) }{ \Ga ( n + \sum_{\ell = 0}^k \vse_\ell ) } \prod_{\ell
= 0}^k \f{ \Ga (x_\ell + \vse_\ell) }{ \Ga (\vse_\ell) }, \qu
\tx{with} \; \vse_0 = \al_0 \; \tx{and} \; 0 < \vse_\ell \leq
\al_\ell, \qu \ell = 1, \cd, k
\]
where $\bs{\vse} = [\vse_0, \vse_1, \cd, \vse_k]^\top$.   Clearly,
\[
\f{ f(\bs{x}) } { g(\bs{x}, \bs{\vse}) } = \f{ \Ga (\sum_{\ell =
0}^k \al_\ell) \; \Ga ( n + \sum_{\ell = 0}^k \vse_\ell ) }{ \Ga
(\sum_{\ell = 0}^k \vse_\ell) \; \Ga ( n + \sum_{\ell = 0}^k
\al_\ell )  }  \prod_{\ell = 1}^k \f{ \Ga (x_\ell + \al_\ell) \; \Ga
(\vse_\ell) }{ \Ga (x_\ell + \vse_\ell) \; \Ga (\al_\ell) }.
\]
For simplicity of notations, define
\[
L (\bs{x}, \bs{\vse}) = \f{ f(\bs{x}) } { g(\bs{x}, \bs{\vse}) }.
\]
Let $\bs{y} = [y_0, y_1, \cd, y_k]^\top$ be a vector such that $y_i
= x_i + 1$ for some $i \in \{ 1, \cd, k \}$ and that $y_\ell =
x_\ell$ for all $\ell \in \{ 1, \cd, k \}$ except $\ell = i$. Then,
\[ \f{  L (\bs{y}, \bs{\vse} )  }{ L (\bs{x}, \bs{\vse} )  } = \f{
x_i + \al_i  }{  x_i + \vse_i  } \geq 1.
\]
Making use of this observation and by an inductive argument, we have
that for $\bs{z} = [z_0, z_1, \cd, z_k]^\top$ such that $x_\ell \leq
z_\ell$ for $\ell = 1, \cd, k$, it must be true that
\[
\f{ L (\bs{z}, \bs{\vse} )  }{ L (\bs{x}, \bs{\vse} )  } \geq 1.
\]
It follows that
\[
\f{ f(\bs{x}) } { g(\bs{x}, \bs{\vse}) } \leq \Lm ( \bs{\vse} ) \qu \fa \bs{\vse} \in \bs{ \varTheta  } \; \tx{provided that} \; x_\ell \leq
z_\ell, \; \ell = 1, \cd, k,
\]
where
\[
\Lm ( \bs{\vse} ) = \f{ \Ga (\sum_{\ell = 0}^k \al_\ell) \; \Ga ( n + \sum_{\ell = 0}^k \vse_\ell ) }{ \Ga (\sum_{\ell = 0}^k \vse_\ell) \; \Ga
( n + \sum_{\ell = 0}^k \al_\ell )  }  \prod_{\ell = 1}^k \f{ \Ga (z_\ell + \al_\ell) \; \Ga (\vse_\ell) }{ \Ga (z_\ell + \vse_\ell) \; \Ga
(\al_\ell) },
\]
$\bs{ \varTheta  }$ is the set of vectors $\bs{\vse} = [\vse_0, \vse_1, \cd, \vse_k]^\top$ such that $\vse_0 = \al_0$ and $0 < \vse_\ell \leq
\al_\ell, \; \ell = 1, \cd, k$.  This implies that
\[
f ( \bs{\mcal{X}} ) \; \bb{ I }_{  \{ \bs{\mcal{X}} \preccurlyeq \bs{z} \} } \leq \Lm ( \bs{\vse} ) \; g ( \bs{\mcal{X} } , \bs{\vse} )  \qu \fa
\bs{\vse} \in \bs{ \varTheta  },
\]
where $\bs{\mcal{X}} = [X_0, X_1, \cd, X_k]^\top$ and $\bs{\mcal{X}} \preccurlyeq \bs{z}$ means $X_\ell \leq z_\ell, \; \ell = 1, \cd, k$.  By
virtue of Theorem \ref{ThM888}, we have $\Pr \li \{ X_\ell \leq z_\ell, \; \ell = 1, \cd, k  \ri \} = \Pr \{  \bs{\mcal{X}} \preccurlyeq \bs{z}
\} \leq \inf_{ \bs{\vse} \in \bs{ \varTheta  } } \Lm ( \bs{\vse} )$.

As a consequence of the assumption that $0 < z_\ell \leq \f{n
\al_\ell}{\sum_{i = 0}^k \al_i}$ for $\ell = 1, \cd, k$, we have
that \[ \se_\ell = \f{\al_0 z_\ell} {  n - \sum_{i=1}^k z_i } \leq
\f{\al_0 \f{n \al_\ell}{\sum_{i = 0}^k \al_i}} {  n - \sum_{\ell
=1}^k \f{n \al_\ell}{\sum_{i = 0}^k \al_i} } =  \al_\ell
\]
for $\ell = 1, \cd, k$.  Define $\bs{\se} = [\se_0, \se_1, \cd,
\se_k]^\top$.  Then, $\bs{\se} \in \bs{ \varTheta  }$ and
\[
\Pr \li \{ X_\ell \leq z_\ell, \; \ell = 1, \cd, k  \ri \}  \leq
\inf_{ \bs{\vse} \in \bs{ \varTheta  } } \Lm ( \bs{\vse} ) \leq \Lm
( \bs{\se} ).
\]
This completes the proof of the theorem.

\subsection{Proof of Theorem \ref{Wishart}} \la{Wishartapp}

To apply the LR method to show inequality (\ref{invGaM}), we introduce a family of probability density functions
\[
g (\bs{x}, \bs{\vse}) = \f{| \bs{\vse} |^{\al}}{ \ba^{p \al}   \Ga_p
(\al ) } | \bs{x} |^{- \al -( p + 1) \sh 2} \exp \li (  - \f{1}{\ba}
\tx{tr} ( \bs{\vse} \bs{x}^{-1} )\ri ),
\]
where $\bs{\vse}$ is a positive-definite real matrix of size $p
\times p$ such that $\bs{\vse} \preccurlyeq \bs{\Psi}$.  Note that
\[
\f{ f(\bs{x}) }{ g(\bs{x}, \bs{\vse}) } = \f{ | \bs{\Psi} |^{\al} }{
| \bs{\vse} |^{\al} } \exp \li (  - \f{1}{\ba} \tx{tr} ( [ \bs{\Psi}
- \bs{\vse}]  \bs{x}^{-1} )\ri ).
\]
For positive definite matrices $\bs{x}$ and $\bs{z}$ such that
$\bs{x} \preccurlyeq \bs{z}$, we have
\[
\tx{tr} ( [ \bs{\Psi} - \bs{\vse}]  \bs{x}^{-1} ) \geq \tx{tr} ( [
\bs{\Psi} - \bs{\vse}]  \bs{z}^{-1} )
\]
as a consequence of $\bs{\vse} \preccurlyeq \bs{\Psi}$.  If follows
that
\[
\f{ f(\bs{x}) }{ g(\bs{x}, \bs{\vse}) } \leq \Lm (\bs{\vse})
\]
for $\bs{\vse} \preccurlyeq \bs{\Psi}$ and $\bs{x} \preccurlyeq
\bs{z}$, where \[ \Lm (\bs{\vse}) = \f{ | \bs{\Psi} |^{\al} }{ |
\bs{\vse} |^{\al} } \exp \li (  - \f{1}{\ba} \tx{tr} ( [ \bs{\Psi} -
\bs{\vse}]  \bs{z}^{-1} )\ri ).
\]
Hence,
\[
f(\bs{X}) \; \bb{I}_{ \{  \bs{X} \preccurlyeq \bs{z}  \}  }  \leq
\Lm (\bs{\vse}) \; g (\bs{X}, \bs{\vse} ) \qu \tx{provided that
$\bs{\vse} \preccurlyeq \bs{\Psi}$}.
\]
 By virtue of Theorem \ref{ThM888}, we have $\Pr \{  \bs{X} \preccurlyeq \bs{z}
 \} \leq \inf_{ \bs{\vse} \preccurlyeq \bs{\Psi} }  \Lm
 (\bs{\vse})$.  In particular, taking $\bs{z} = \ro \bb{E} [ \bs{X}
 ] = \f{2\ro}{\ba}  \f{ \bs{\Psi} }{ 2 \al - p - 1 }$ and $\bs{\vse} = \f{\ba}{2} (2
\al - p - 1) \bs{z}$, we have \bee
 \Pr \li \{ \bs{X} \preccurlyeq \ro \bs{\Up} \ri \} & = &  \Pr \li \{ \bs{X} \preccurlyeq
 \bs{z}  \ri \}\\
 & \leq & \f{ | \bs{\Psi} |^{\al} }{  | \f{\ba}{2} (2 \al - p - 1)
\bs{z} |^{\al} } \exp \li (  -
\f{1}{\ba} \tx{tr} ( [ \bs{\Psi} - \f{\ba}{2} (2 \al - p - 1) \bs{z} ] \bs{z}^{-1} ) \ri )\\
& = &  \f{   \exp (  \f{p}{2} ( 2 \al - p - 1 )   ) | \bs{\Psi}
|^{\al}  }{  [ \f{\ba}{2} (2 \al - p - 1)  ]^{p \al}  | \bs{z}
|^{\al} } \exp \li ( - \f{1}{\ba} \tx{tr} ( \bs{\Psi}  \bs{z}^{-1} )
\ri ) \\
& = & \f{ 1  }{ \ro^{p \al} } \exp \li ( - \f{p}{2} \li ( \f{1}{\ro}
- 1 \ri ) ( 2 \al - p - 1 ) \ri ).
 \eee
This completes the proof of the theorem.

\subsection{Proof of Theorem \ref{mulnormal}} \la{mulnormalapp}

For simplicity of notations, let
\[
\bs{\mcal{X}} = [\bs{X}_1, \cd, \bs{X}_n].
\]
Let
\[
\mscr{X} = [ \bs{x}_1, \cd, \bs{x}_n],
\]
where $\bs{x}_1, \cd, \bs{x}_n$ are vectors of dimension $k$.  Since
$\bs{X}_1, \cd, \bs{X}_n$ are identical and independent, the joint
probability density of $\bs{\mcal{X}}$ is
\[
f_{ \bs{\mcal{X}} } ( \mscr{X} ) = \prod_{i = 1}^n \li [( 2 \pi)^{-
k \sh 2} | \bs{\Si} |^{-1 \sh 2} \exp \li ( - \f{1}{2} (\bs{x}_i -
\bs{\mu})^{\top} \bs{\Si}^{-1} (\bs{x}_i - \bs{\mu})  \ri ) \ri ]^n.
\]
To apply the LR method to show (\ref{multinormal}), we introduce a family of probability density functions
\[
g_{ \bs{\mcal{X}} } ( \mscr{X}, \bs{\vse} ) = \prod_{i = 1}^n \li [(
2 \pi)^{- k \sh 2} | \bs{\Si} |^{-1 \sh 2} \exp \li ( - \f{1}{2}
(\bs{x}_i - \bs{\vse})^{\top} \bs{\Si}^{-1} (\bs{x}_i - \bs{\vse})
\ri ) \ri ]^n,
\]
where $\bs{\vse}$ is a vector of dimension $k$ such that
$\bs{\Si}^{-1}  \bs{\vse} \succcurlyeq \bs{\Si}^{-1} \bs{\mu}$. It
can be checked that \bee  \f{ f_{ \bs{\mcal{X}} } ( \mscr{X} ) } {
g_{ \bs{\mcal{X}} } ( \mscr{X}, \bs{\vse} )  } & = & \prod_{i = 1}^n
\exp \li ( (\bs{\mu}^{\top} - \bs{\vse}^{\top}) \bs{\Si}^{-1}
\bs{x}_i + \f{1}{2} [ \bs{\vse}^{\top} \bs{\Si}^{-1} \bs{\vse} -
\bs{\mu}^{\top} \bs{\Si}^{-1} \bs{\mu} ]  \ri )\\
& = & \li [ \exp \li ( (\bs{\mu}^{\top} - \bs{\vse}^{\top})
\bs{\Si}^{-1} \ovl{\bs{x}}_n + \f{1}{2} [ \bs{\vse}^{\top}
\bs{\Si}^{-1} \bs{\vse} - \bs{\mu}^{\top} \bs{\Si}^{-1} \bs{\mu} ]
\ri ) \ri ]^n,  \eee where
\[
\ovl{\bs{x}}_n = \f{ \sum_{i=1}^n \bs{x}_i }{n}.
\]
As a consequence of $\bs{\Si}^{-1}  \bs{\vse} \succcurlyeq
\bs{\Si}^{-1} \bs{\mu}$, we have that
\[
(\bs{\mu}^{\top} - \bs{\vse}^{\top}) \bs{\Si}^{-1} \bs{u}\leq
(\bs{\mu}^{\top} - \bs{\vse}^{\top}) \bs{\Si}^{-1} \bs{v}
\]
for arbitrary vectors $\bs{u}$ and $\bs{v}$ such that $\bs{v}
\succcurlyeq \bs{u}$.  This implies that for $\bs{\vse}$ such that
$\bs{\Si}^{-1}  \bs{\vse} \succcurlyeq \bs{\Si}^{-1} \bs{\mu}$,
\[
\f{ f_{ \bs{\mcal{X}} } ( \mscr{X} ) } { g_{ \bs{\mcal{X}} } (
\mscr{X}, \bs{\vse} )  }  \leq \Lm ( \bs{\vse} ) \qu \tx{ provided
that $\ovl{\bs{x}}_n \succcurlyeq \bs{z}$,  }
\]
where \[ \Lm ( \bs{\vse} ) = \li [ \exp \li ( (\bs{\mu}^{\top} -
\bs{\vse}^{\top}) \bs{\Si}^{-1} \bs{z} + \f{1}{2} [ \bs{\vse}^{\top}
\bs{\Si}^{-1} \bs{\vse} - \bs{\mu}^{\top} \bs{\Si}^{-1} \bs{\mu} ]
\ri ) \ri ]^n.
\]
It follows that
\[
f_{ \bs{\mcal{X}} } ( \bs{\mcal{X}} ) \; \bb{ I }_{  \{ \ovl{\bs{X}}_n \succcurlyeq \bs{z} \} } \leq \Lm ( \bs{\vse} ) \; g_{ \bs{\mcal{X}} }  (
\bs{\mcal{X} } , \bs{\vse} )
\]
for any $\bs{\vse}$ such that $\bs{\Si}^{-1}  \bs{\vse} \succcurlyeq \bs{\Si}^{-1} \bs{\mu}$.  By virtue of By virtue of Theorem \ref{ThM888},
we have  \bee \Pr \{ \ovl{\bs{X}}_n \geq \bs{z}  \} \leq \inf_{ \bs{\Si}^{-1}  \bs{\vse} \succcurlyeq \bs{\Si}^{-1} \bs{\mu} } \Lm ( \bs{\vse}
). \eee By differentiation, it can be shown that the infimum is attained at $ \bs{\vse} = \bs{z}$. Hence, $\Pr \{ \ovl{\bs{X}}_n \geq \bs{z}  \}
\leq \Lm (\bs{\vse})$. This completes the proof of the theorem.

\subsection{Proof of Theorem \ref{multiPareto}} \la{multiParetoapp}

Let $\bs{\mcal{X}}$ denote the random matrix of size $k \times n$
such that the $j$-th column is $\bs{\fra{X}}_j$. Let $\bs{x}$ denote
a matrix of size $k \times n$ such that the $j$-th column is
$[x_{1j}, x_{2j}, \cd, x_{kj}]^\top$.   Then, the joint probability
density function of $\bs{\mcal{X}}$  is
\[
f_{ \bs{\mcal{X}} } ( \bs{x} ) = \prod_{j=1}^n \li [ \li ( \prod_{i
= 1}^k \f{\al + i - 1}{\ba_i} \ri )  \li ( 1 - k + \sum_{i=1}^k
\f{x_{ij}}{\ba_i} \ri )^{- (\al + k) } \ri ].
\]
To apply the LR method, we introduce a family of probability density
functions
\[
g_{ \bs{\mcal{X}} } ( \bs{x}, \vse ) = \prod_{j=1}^n \li [ \li (
\prod_{i = 1}^k \f{\vse + i - 1}{\ba_i} \ri )  \li ( 1 - k +
\sum_{i=1}^k \f{x_{ij}}{\ba_i} \ri )^{- (\vse + k) } \ri ], \qu \vse
> \al.
\]
Clearly,
\[
\f{f_{\bs{\mcal{X}}} (\bs{x})} { g_{\bs{\mcal{X}}} (\bs{x}, \vse)} =
\li ( \prod_{i = 1}^k \f{\al + i - 1}{\vse + i - 1} \ri )^n \li [
\prod_{j=1}^n  \li ( 1 - k + \sum_{i=1}^k \f{x_{ij}}{\ba_i} \ri )
\ri ]^{\vse - \al}.
\]
Using the fact that the geometric mean does not exceed the
arithmetic mean, we have that
\[
\prod_{j=1}^n  \li ( 1 - k + \sum_{i=1}^k \f{x_{ij}}{\ba_i} \ri )
\leq \li ( 1 - k + \sum_{i=1}^k \f{u_i}{\ba_i} \ri )^n,
\]
where
\[
u_i = \f{1}{n} \sum_{j=1}^n x_{ij}.
\]
Hence,
\[
\f{f_{\bs{\mcal{X}}} (\bs{x})} { g_{\bs{\mcal{X}}} (\bs{x}, \vse)}
\leq \li [  \li ( \prod_{i = 1}^k \f{\al + i - 1}{\vse + i - 1} \ri
) \li ( 1 - k + \sum_{i=1}^k \f{u_i}{\ba_i} \ri )^{\vse - \al} \ri
]^n.
\]
This implies that
\[
f_{\bs{\mcal{X}}} ( \bs{\mcal{X}}) \; \bb{I}_{\{  \ovl{\bs{\fra{X}}
}_n \preceq \bs{z} \}} \leq \Lm (\vse) \; g_{\bs{\mcal{X}}} (
\bs{\mcal{X}}, \vse), \qu \fa \vse > \al,
\]
where
\[
\Lm (\vse) = \li [  \li ( \prod_{i = 1}^k \f{\al + i - 1}{\vse + i -
1} \ri ) \li ( 1 - k + \sum_{i=1}^k \f{z_i}{\ba_i} \ri )^{\vse -
\al} \ri ]^n.
\]
By virtue of Theorem \ref{ThM888}, we have  $\Pr \{
\ovl{\bs{\fra{X}} }_n \preceq \bs{z} \} \leq \Lm (\se)$ for any $\se
> \al$. This proves the first statement.

Note that if (\ref{best9988}) holds, then $\se > \al$ for $\se$
satisfying (\ref{bestp88}). Moreover, by differentiation, it can be
shown that $\Pr \{ \ovl{\bs{\fra{X}} }_n \preceq \bs{z} \} \leq
\inf_{\vse > \al} \Lm (\vse) = \Lm (\se)$. This proves statement
(II).

For $\se$ satisfying (\ref{mean889}), it must be true that $\se >
\al$ as a consequence of the assumption that $\al > 1$ and $\f{1}{k}
\sum_{i=1}^k \f{z_i}{\ba_i} < \f{\al}{\al - 1}$. This proves
statement (III). The proof of the theorem is thus completed.

\end{document}